\documentclass[a4paper,UKenglish,cleveref]{lipics-v2019-hacked}

\usepackage{thmtools}
\usepackage{thm-restate}
\usepackage[utf8]{inputenc}
\usepackage{amsmath,amssymb,amsthm,xspace,xcolor,graphicx}
\usepackage{changepage}
\usepackage{enumitem}
\usepackage{tikz}
\usetikzlibrary{fit,automata,arrows}
\usetikzlibrary{patterns,hobby}
\usepackage{caption}
\usepackage[labelformat=simple]{subcaption}

\usepackage{float}

\title{Comparing Width Parameters on Graph Classes}

\titlerunning{Comparing Width Parameters on Graph Classes}

\author{Nick Brettell}{School of Mathematics and Statistics, Victoria University of Wellington, New Zealand}{nick.brettell@vuw.ac.nz}{https://orcid.org/0000-0002-1136-418X}{}
\author{Andrea Munaro}{Department of Mathematical, Physical and Computer Sciences, University of Parma, Italy}{andrea.munaro@unipr.it}{https://orcid.org/0000-0003-1509-8832}{}
\author{Dani\"el Paulusma}{Department of Computer Science, Durham University, UK}{daniel.paulusma@durham.ac.uk}{https://orcid.org/0000-0001-5945-9287}{}
\author{Shizhou Yang}{School of Mathematics and Physics, Queen's University Belfast, UK}{syang22@qub.ac.uk}{}{}

\authorrunning{N. Brettell, A. Munaro, D. Paulusma and S. Yang}


\acknowledgements{The first author was supported by a New Zealand Marsden fund. The second author is a member of the Gruppo Nazionale Calcolo Scientifico-Istituto Nazionale di Alta Matematica (GNCS-INdAM).}

\nolinenumbers 

\hideLIPIcs

\newtheorem{open}{Open Problem}
\newtheorem{conj}[theorem]{Conjecture}

\theoremstyle{definition}
\newtheorem{observation}[theorem]{Observation}
\newtheorem{subobservation}{Observation}[theorem]
\newtheorem{subclaim}[subobservation]{Claim}

\newcommand{\cutmim}{\mathrm{cutmim}}
\newcommand{\cutsim}{\mathrm{cutsim}}
\newcommand{\mmw}{\mathrm{mmw}}
\newcommand{\mimw}{\mathrm{mimw}}
\newcommand{\simw}{\mathrm{simw}}
\newcommand{\cw}{\mathrm{cw}}
\newcommand{\tw}{\mathrm{tw}}
\newcommand{\tww}{\mathrm{tww}}
\newcommand{\nlcw}{\mathrm{nlcw}}
\newcommand{\bw}{\mathrm{bw}}
\newcommand{\tin}{\mathrm{tree}\textnormal{-}\alpha}

\renewcommand\cref[1]{\Cref{#1}}

\DeclareMathOperator{\rw}{rw}

\begin{document}

\maketitle

\begin{abstract}
We study how the relationship between non-equivalent width parameters changes once we restrict to some special graph class. As width parameters we consider treewidth, clique-width, twin-width, mim-width, sim-width and tree-independence number, whereas as graph classes we consider $K_{t,t}$-subgraph-free graphs, line graphs and their common superclass, for $t \geq 3$, of $K_{t,t}$-free graphs. For $K_{t,t}$-subgraph-free graphs, we extend a known result of Gurski and Wanke (2000) and provide a complete comparison, showing in particular that treewidth, clique-width, mim-width, sim-width and tree-independence number are all equivalent. For line graphs, we extend a result of Gurski and Wanke (2007) and also provide a complete comparison, showing in particular that clique-width, mim-width, sim-width and tree-independence number are all equivalent, and bounded if and only if the class of root graphs has bounded treewidth.  For $K_{t,t}$-free graphs, we provide an almost-complete comparison, leaving open only one missing case. We show in particular that $K_{t,t}$-free graphs of bounded mim-width have bounded tree-independence number, and obtain structural and algorithmic consequences of this result, such as a proof of a special case of a recent conjecture of Dallard, Milani\v{c} and \v{S}torgel. Finally, we consider the question of whether boundedness of a certain width parameter is preserved under graph powers. We show that this question has a positive answer for sim-width precisely in the case of odd powers.
\end{abstract} 

\section{Introduction}\label{sec:intro}

Width parameters play an important role both in structural and algorithmic graph theory. 
Two width parameters may differ in strength. To explain this, we say that a width parameter $p$ {\it dominates} a width parameter~$q$ if there exists a function~$f$ such that $p(G)\leq f(q(G))$ for every graph~$G$.
If $p$ dominates $q$ but $q$ does not dominate $p$, then $p$ is {\it more powerful} than $q$.
If $p$ dominates $q$ and $q$ dominates $p$, then $p$ and $q$ are {\it equivalent}. In particular, if two equivalent parameters $p$ and $q$ admit linear functions witnessing this, we say that $p$ and $q$ are \textit{linearly equivalent} (in other words, one is a constant factor approximation of the other). If neither $p$ dominates $q$ nor $q$ dominates $p$, then $p$ and $q$ are {\it incomparable}. A width parameter~$p$ is {\it bounded} on a graph class ${\cal G}$ if there exists a constant~$c$ such that $p(G)\leq c$ for every $G\in {\cal G}$. Note that if $p$ and $q$ are two equivalent width parameters then, for every graph class~${\cal G}$, the parameter $p$ is bounded on~${\cal G}$ if and only if $q$ is bounded on ${\cal G}$. 

We can define all the above notions with respect to special graph classes (instead of the class of all graphs) analogously. In particular, two width parameters $p$ and $q$ are {\it equivalent for some graph class} ${\cal G}$ if, for every subclass ${\cal G}'$ of~${\cal G}$, the parameter $p$ is bounded on ${\cal G}'$ if and only if $q$ is bounded on ${\cal G}'$. This definition leads to a natural research question:

\medskip
\noindent
{\it For which graph classes does the relationship between two non-equivalent width parameters $p$ and $q$ change?}

\medskip
\noindent
For example, two width parameters $p$ and $q$ might be incomparable in general, but when restricted to some special graph class, one of them could dominate the other. Or $p$ might be more powerful than $q$ in general, but when restricted to some special graph class, $p$ and $q$ could become equivalent.

In order to explain our results and how they embed in the literature, we first need to introduce some terminology. For a graph $H$, a graph $G$ is {\it $H$-subgraph-free} if $G$ cannot be modified into $H$ by a sequence of edge deletions and vertex deletions, whereas $G$ is {\it $H$-free} if $G$ cannot be modified into $H$ by a sequence of vertex deletions. For a set of graphs $\{H_1,\ldots,H_k\}$ for some $k\geq 1$, a graph $G$ is \emph{$(H_1,\ldots,H_k)$-free} (or \emph{$(H_1,\ldots,H_k)$-subgraph-free}) if $G$ is $H_i$-free (or $H_i$-subgraph-free) for every $i\in \{1,\ldots,k\}$. For an integer $s$, we let $K_s$ denote the complete graph on $s$ vertices. For integers $s$ and $t$, we let $K_{s,t}$ denote the complete bipartite graph whose partition classes have size~$s$ and $t$. We also say that a complete bipartite graph is a \emph{biclique}. The \emph{line graph} $L(G)$ of a graph~$G$ has vertex set~$E$ and an edge between two vertices~$e_1$ and~$e_2$ if and only if~$e_1$ and~$e_2$ have a common endpoint in~$G$.

\begin{figure}
\begin{center}
\includegraphics[scale=0.8]{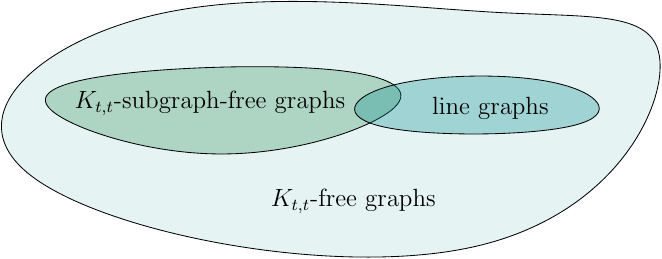}
\caption{Inclusion diagram of the three main graph classes considered in the paper, where $t \geq 3$.}
\label{f-classes}
\end{center}
\end{figure}

Our starting point consists of two well-known results of Gurski and Wanke characterizing clique-width on $K_{t,t}$-subgraph-free graphs \cite{GW00} and line graphs \cite{GW07} in terms of treewidth. The first result \cite{GW00} asserts that
clique-width and treewidth are equivalent for the class of $K_{t,t}$-subgraph-free graphs, for every $t \geq 2$.
This implies that rank-width\footnote{A parameter equivalent to clique-width for the class of all graphs (see \Cref{s-pre}).} and treewidth are also equivalent for this graph class, and an improved upper bound for treewidth which is polynomial in rank-width was obtained by Fomin, Oum and Thilikos~\cite{FOT10}. The second result of Gurski and Wanke~\cite{GW07} asserts that, for any graph class $\mathcal{G}$, the class of line graphs of graphs in $\mathcal{G}$ has bounded clique-width if and only if $\mathcal{G}$ has bounded treewidth. In particular, although clique-width and treewidth are equivalent for the class of $K_{t,t}$-subgraph-free graphs (one of the least restrictive classes of sparse graphs), clique-width is more powerful than treewidth for line graphs (one of the most prominent classes of dense graphs).

In this paper we answer the question above (except for one open case) for the following three graph classes and six width parameters: 

\begin{itemize}
\item graph classes (see \Cref{f-classes}): $K_{t,t}$-subgraph-free graphs, line graphs, and their common proper superclass (for $t \geq 3$) of $K_{t,t}$-free graphs. 
\item width parameters (see \Cref{f-width}): treewidth ($\tw$), clique-width ($\cw$), twin-width ($\tww$), mim-width ($\mimw$), sim-width ($\simw$) and tree-independence number ($\tin$).    
\end{itemize}

\noindent
Apart from the fact that $K_{t,t}$-free graphs (with $t\geq 3$) generalize $K_{t,t}$-subgraph-free graphs and line graphs, our motivation for investigating $K_{t,t}$-free graphs is that large bicliques are obstructions to small tree-independence number. Our main results are summarized in \Cref{figkttsubgraphfree,f-line,kttrelfig} and will be explained in detail in \Cref{s-results}.

We also consider an extensively studied question concerning the behaviour of width parameters under graph powers. Given a graph $G$, the \textit{$r$-th power} of $G$ is the graph $G^r$ whose vertex set is that of $G$ and where two distinct vertices are adjacent if and only if they are at distance at most $r$ in $G$. 

\medskip
\noindent
{\it Given a width parameter $p$, is boundedness of $p$ preserved under graph powers?}

\medskip
\noindent
The behaviour of each parameter $p \in \{\tw, \cw, \tww, \mimw, \tin\}$ is known \cite{BKTW21,JKST19,LMM24,ST07}, and we show that the question has a positive answer for $p = \simw$ precisely in the case of odd powers. 

Before explaining our results in detail and analyzing their consequences, we recall how the six parameters we consider are related to each other on the class of all graphs (we refer the reader to \Cref{s-pre} for the definitions of these parameters and for the computable functions showing that one dominates another). Kang et al.~\cite{KKST17} proved that sim-width is more powerful than mim-width, while Vatshelle~\cite{Vat12} showed that mim-width is more powerful than clique-width. Courcelle and Olariu~\cite{CO00} showed that clique-width is more powerful than treewidth, whereas
Bonnet et al.~\cite{BKTW21} proved that twin-width is more powerful than clique-width. Munaro and Yang~\cite{MY23} showed that sim-width is more powerful than tree-independence number, which in turn is more powerful than treewidth, as shown independently by Yolov~\cite{Yol18} and Dallard, Milani\v{c} and \v{S}torgel~\cite{DMS24}.
Regarding incomparability, it is known that complete bipartite graphs have clique-width~$2$, but unbounded tree-independence number~\cite{DMS24}. Moreover, chordal graphs have tree-independence number~$1$~\cite{DMS24} but unbounded twin-width~\cite{BGKTW21} and unbounded mim-width~\cite{KKST17}. 
In fact, even interval graphs have unbounded twin-width~\cite{BGKTW21}, but in contrast to their superclass of chordal graphs they have bounded mim-width~\cite{BV13}.
Finally, walls have bounded twin-width \cite{BKTW21}, but using a result of \cite{KKST17} it can be easily shown that they have unbounded sim-width (see Section~\ref{s-pre} for a proof). 
We conclude that twin-width is incomparable to mim-width, sim-width and tree-independence number, and that the latter parameter is incomparable to mim-width and clique-width.
%
%
%
\Cref{f-width} gives a full picture of all relationships between the six parameters we just described.

\begin{figure}[htbp]
\centering
\begin{subfigure}[b]{\textwidth}
\captionsetup{justification=justified}
\begin{tikzpicture}[yscale=.5]
\node at (4.5,9.8) (tin) {tree-independence number};
\node at (-3.5,11) (twin) {twin-width};
\node at (2.58,7.65) (ref6) {\cite{DMS24,Yol18}};
\node at (-2.75,9) (ref4) {\cite{BKTW21}};
\node at (1.6,11.05) (ref7) {\cite{MY23}};
\node at (0,11) (p) {sim-width};
\node at (-0.4,10.1) (ref1) {\cite{KKST17}};
\node at (-3,8.5) (b) {};
\node at (0,9) (wc) {mim-width};
\node at (-0.4,8.1) (ref2) {\cite{Vat12}};
\node at (-6,7) (cb) {};
\node at (0,7) (c) {clique-width};
\node at (-0.4,6.1) (ref3) {\cite{CO00}};
\node at (3,7) (dh) {};
\node at (-6,5.5) (bp) {};
\node at (-1.5,5.5) (s) {};
\node at (0,5) (i) {treewidth};
\draw [->, line width=0.3mm](p) -- (tin.175);
\draw [->, line width=0.3mm](twin.330) -- (c.173);
\draw [->, line width=0.3mm](tin.188) -- (i.18);
\draw [->, line width=0.3mm] (p) -- (wc);
\draw [->, line width=0.3mm] (wc) -- (c);
\draw [->, line width=0.3mm] (c) -- (i);
\end{tikzpicture}
\caption{The relationships between the six different width parameters on the class of all graphs.}\label{f-width}
\end{subfigure}

\vspace{0.5cm}

\begin{subfigure}[b]{\textwidth}
\captionsetup{justification=justified}
\begin{tikzpicture}[yscale=.45]
\node at (4.5,9.8) (tin) {tree-independence number};
\node at (0,11) (p) {sim-width};
\node at (-3,8.5) (b) {};
\node at (0,9) (wc) {mim-width};
\node at (-6,7) (cb) {};
\node at (0,7) (c) {clique-width};
\node at (3,7) (dh) {};
\node at (-6,5.5) (bp) {};
\node at (-1.5,5.5) (s) {};
\node at (0,5) (i) {treewidth};
\node at (-3.5,11) (twin) {twin-width};
\draw [->, line width=0.3mm](twin) -- (p);
\draw [<->, line width=0.3mm](tin.175) -- (p);
\draw [<->, line width=0.3mm] (p) -- (wc);
\draw [<->,  line width=0.3mm] (wc) -- (c);
\draw [<->, line width=0.3mm] (c) -- (i);
\draw [<->, line width=0.3mm](tin.188) -- (i.18);
\end{tikzpicture}
\caption{The relationships between the six different width parameters when restricted to $(K_s,K_{t,t})$-free graphs for some $s\geq 3$ and $t\geq 2$, or to $K_{t,t}$-subgraph-free graphs for some $t\geq 2$. See \Cref{t-main} and \Cref{c-main}, respectively.}\label{figkttsubgraphfree}
\end{subfigure}

\vspace{0.5cm}

\begin{subfigure}[b]{\textwidth}
\captionsetup{justification=justified}
\begin{tikzpicture}[yscale=.5]
\node at (4.5,9.8) (tin) {tree-independence number};
\node at (0,11) (p) {sim-width};
\node at (-3,8.5) (b) {};
\node at (0,9) (wc) {mim-width};
\node at (-6,7) (cb) {};
\node at (0,7) (c) {clique-width};
\node at (3,7) (dh) {};
\node at (-6,5.5) (bp) {};
\node at (-1.5,5.5) (s) {};
\node at (0,5) (i) {treewidth};
\node at (-3.5,11) (twin) {twin-width};
\draw [->, line width=0.3mm](twin) -- (p);
\draw [<->, line width=0.3mm](tin.175) -- (p);
\draw [<->, line width=0.3mm] (p) -- (wc);
\draw [<->, line width=0.3mm] (wc) -- (c);
\draw [<->, line width=0.3mm ](tin.188) -- (c.18);
\draw [->, line width=0.3mm] (c) -- (i);
\end{tikzpicture}
\caption{The relationships between the six different width parameters when restricted to line graphs. See \Cref{t-main3}.}
\label{f-line}
\end{subfigure}

\vspace{0.5cm}

\begin{subfigure}[b]{\textwidth}
\captionsetup{justification=justified}
\begin{tikzpicture}[yscale=.5]
  \node at (4.5,9.8) (tin) {tree-independence number};
\node at (0,11) (p) {sim-width};
\node at (-3,8.5) (b) {};
\node at (0,9) (wc) {mim-width};
\node at (-6,7) (cb) {};
\node at (0,7) (c) {clique-width};
\node at (3,7) (dh) {};
\node at (-6,5.5) (bp) {};
\node at (-1.5,5.5) (s) {};
\node at (0,5) (i) {treewidth};
\node at (-3.5,11) (twin) {twin-width};
\draw [->, line width=0.3mm](p) -- (tin.175);
\path (p) edge[bend left,red,<-, line width=0.3mm] node [above] {?} (tin);
\draw [->, line width=0.3mm](twin.330) -- (c.173);
\draw [->, line width=0.3mm](tin) -- (wc);
\draw [->, line width=0.3mm] (p) -- (wc);
\draw [->, line width=0.3mm] (wc) -- (c);
\draw [->, line width=0.3mm] (c) -- (i);
\end{tikzpicture}
\caption{The relationships between the six different width parameters when restricted to $K_{t,t}$-free graphs for some $t \ge 2$. See \Cref{t-main2}. The red arrow illustrates the one remaining open case.}
\label{kttrelfig}
\end{subfigure}

\caption{The relationships between the six different width parameters considered in the paper when restricted to different graph classes. A directed arrow from parameter $p$ to parameter $q$ indicates that $p$ dominates $q$, whereas a bidirected arrow indicates that $p$ and $q$ are equivalent. Although not explicitly stated in \Cref{t-main,t-main2,t-main3}, all functions showing that a certain parameter $p$ dominates another parameter $q$ for a certain graph class~$\mathcal{G}$ are computable and, as we will show in particular for the new results, can be obtained from the corresponding proofs.}\label{fig:gadham}
\end{figure}

\subsection{Our Results}\label{s-results}

In the following four sections we highlight our main results. We then observe some of their algorithmic and structural consequences and connections with known results.

\subsubsection{$\mathbf{K_{t,t}}$-Subgraph-Free Graphs} 

The class of $K_{t,t}$-subgraph-free graphs contains well-known sparse graph classes: 
for example, every degenerate graph class and every nowhere dense graph class is a (proper) subclass of the class of $K_{t,t}$-subgraph-free graphs for some $t\geq 3$ (see~\cite{TV19} for a proof).
Gurski and Wanke~\cite{GW00} proved that for every $t\geq 2$, clique-width and treewidth become equivalent for the class of $K_{t,t}$-subgraph-free graphs. In \Cref{s-main} we generalize and extend their result as follows.

\begin{restatable}{theorem}{tmain}
\label{t-main}
For every $s\geq 3$ and $t\geq 2$, when restricted to $(K_s,K_{t,t})$-free graphs, sim-width, mim-width, clique-width, treewidth and tree-independence number are equivalent, whereas twin-width is more powerful than any of these parameters.
\end{restatable}

The relationships in \cref{t-main} are displayed in Figure~\ref{figkttsubgraphfree}. \Cref{t-main} shows in particular that treewidth, clique-width, mim-width, sim-width and tree-independence number are equivalent for the class of $(K_s,K_{t, t})$-free graphs.
When $s \geq 2t$, the class of $(K_s,K_{t,t})$-free graphs contains the class of $K_{t,t}$-subgraph-free graphs, so these parameters are also equivalent for the class of $K_{t,t}$-subgraph-free graphs. Thus, \cref{t-main} indeed generalizes and extends the result of Gurski and Wanke~\cite{GW00}.

\begin{corollary}\label{c-main} 
For every $t\geq 2$, when restricted to $K_{t,t}$-subgraph-free graphs, sim-width, mim-width, clique-width, treewidth and tree-independence number are equivalent, whereas twin-width is more powerful than any of these parameters.
\end{corollary}

The relationships in \cref{c-main} are also displayed in Figure~\ref{figkttsubgraphfree}. In \Cref{s-main} we also provide a counterexample to a purported result of Vatshelle~\cite[Lemma~4.3.9]{Vat12} stating that a $d$-degenerate bipartite graph with a matching of size $\mu$ has an induced matching of size at least $\mu / (d+1)$. This result was used to obtain a lower bound on the mim-width of grids \cite[Lemma~4.3.10]{Vat12} and later to obtain a lower bound for mim-width that is linear in treewidth for $d$-degenerate graphs \cite{KKST17,Men18}. We explain how to obtain alternative lower bounds in these situations.

\subsubsection{Line Graphs} 

We use \Cref{c-main} in the proof of our next theorem, which concerns line graphs.
For a graph class ${\cal G}$, we let $L(\mathcal{G})$ denote the class of line graphs of graphs in $\mathcal{G}$.
Some years after~\cite{GW00}, Gurski and Wanke~\cite{GW07} proved that a class of graphs~${\cal G}$ has bounded treewidth if and only if $L({\cal G})$ has bounded clique-width. We extend this result by proving the following theorem in \Cref{s-main3}.

\begin{restatable}{theorem}{tmainthree}
\label{t-main3}
For a graph class $\mathcal{G}$, the following statements are equivalent:
\begin{enumerate}
\item The class $\mathcal{G}$ has bounded treewidth;
\item The class $L(\mathcal{G})$ has bounded clique-width;
\item The class $L(\mathcal{G})$ has bounded mim-width;
\item The class $L(\mathcal{G})$ has bounded sim-width;
\item The class $L(\mathcal{G})$ has bounded tree-independence number.
\end{enumerate}
Moreover, when restricted to line graphs, sim-width, mim-width, clique-width and tree-independence number are equivalent; twin-width dominates each of these four parameters; and each of the four parameters in turn dominates treewidth.
\end{restatable}

The relationships in \cref{t-main3} are displayed in Figure~\ref{f-line}. The main technical contribution towards proving \Cref{t-main3} is the following result (\Cref{simline}): There exists a non-decreasing unbounded function $f\colon \mathbb{N} \rightarrow \mathbb{N}$ such that, for every graph $G$, $\simw(L(G)) \geq f(\tw(G))$. 

\Cref{t-main3} immediately implies that, for each parameter $p \in \{\cw,\mimw,\simw,\tin\}$, there is no function $f$ such that $p(L(G)) \leq f(p(G))$ for every graph $G$. The same clearly holds for $p = \tw$ and $p = \tww$ as well \cite{BKTW21}. Moreover, \Cref{t-main3} shows that, for each parameter $p \in \{\cw,\mimw,\simw,\tin\}$, there exist functions $f$ and~$g$ such that $f(\tw(G)) \leq p(L(G)) \leq g(\tw(G))$, for every graph $G$. In fact, the proof shows that we can always choose $g$ as a linear function. However, it is not immediately clear what order of magnitude $f$ should have. For $p = \cw$, Gurski and Wanke~\cite{GW07} showed that a linear function suffices. More precisely, for any graph $G$, they showed that the following holds:
\begin{equation}\label{cwline}
\frac{\tw(G)+1}{4} \leq \cw(L(G)) \leq 2\tw(G) + 2.
\end{equation}
But what about the other cases? In \Cref{s-main3}, we answer this question for $p = \mimw$. In fact, instead of treewidth, we consider its linearly equivalent parameter branch-width ($\bw$) and show that the mim-width of a line graph equals, up to a multiplicative constant, the branch-width of the root graph:

\begin{restatable}{theorem}{mimlinebw}
\label{mimlinebw}
For any graph $G$, $\frac{\bw(G)}{25} \leq \mimw(L(G)) \leq \bw(G)$.
\end{restatable}

We remark that a result similar in spirit to \Cref{mimlinebw} and \Cref{cwline} was obtained by Oum~\cite{Oum09}: for any graph $G$, $\rw(L(G))$ is exactly one of $\bw(G), \bw(G) - 1, \bw(G) - 2$, where rank-width ($\rw$) is a width parameter equivalent to clique-width.

It is in general non-trivial to determine the exact values of $p(L(K_{n,m}))$ and $p(L(K_n))$, for some width parameter $p$. For $n\geq 3$, Lucena~\cite{Luc07} showed that $\tw(L(K_{n,n})) = n^2/2 + n/2 - 1$. This result was extended by Harvey and Wood~\cite{HW18}, who showed that $\tw(L(K_{n,m}))$ has order $nm$. Moreover, Harvey and Wood~\cite{HW15} determined the exact value of $\tw(L(K_{n}))$. We conclude \Cref{s-main3} by showing that $\simw(L(K_{n,m})) = \lceil n/3\rceil$, for any $6 < n \leq m$, and use this result to bound $\simw(L(K_{n}))$. 

\subsubsection{$\mathbf{K_{t,t}}$-Free Graphs} 

In \Cref{s-main2} we show the following result concerning $K_{t,t}$-free graphs.

\begin{restatable}{theorem}{tmaintwo}
\label{t-main2}
For every $t\geq 2$, when restricted to $K_{t,t}$-free graphs,
\begin{itemize}
  \item sim-width dominates tree-independence number, tree-independence number and sim-width are more powerful than mim-width, and twin-width is incomparable with these three parameters;
  \item twin-width and mim-width are more powerful than clique-width; and
  \item clique-width is more powerful than treewidth.
\end{itemize}
\end{restatable}

The relationships in \cref{t-main2} are displayed in Figure~\ref{kttrelfig}. We note from Figures~\ref{f-width} and~\ref{kttrelfig} that, on $K_{t,t}$-free graphs, tree-independence number becomes more powerful than mim-width. Moreover, the equivalences from Figure~\ref{figkttsubgraphfree} no longer hold, except for perhaps one possible relationship. That is, we do not know yet if tree-independence number dominates sim-width when restricted to $K_{t,t}$-free graphs. If so, then these parameters become equivalent when restricted to $K_{t,t}$-free graphs. This is the {\it only} missing case in Figure~\ref{kttrelfig} and the following remains open:

\begin{restatable}{open}{kttfreeopen}
\label{kttfreeopen}
Does tree-independence number dominate sim-width for the class of $K_{t,t}$-free graphs? In other words, is it true that every subclass of $K_{t,t}$-free graphs of bounded sim-width has bounded tree-independence number? 
\end{restatable}

The main ingredient in the proof of \Cref{t-main2} is the following result, which shows in particular that tree-independence number dominates mim-width on $K_{t,t}$-free graphs. Its other consequences will be described in \Cref{conseqtinmim}.

\begin{restatable}{theorem}{tinmim}
\label{tinmim}
Let $n$ and $m$ be positive integers. Let $G$ be a $K_{n,m}$-free graph and let $(T,\delta)$ be a branch decomposition of $G$ with $\mimw_{G}(T, \delta) < k$. Then we can construct a tree decomposition of $G$ with independence number less than $6(2^{n+k-1} + mk^{n+1})$ in $O(|V(G)|^{mk^{n}+4})$ time. In particular, $\tin(G) < 6(2^{n+k-1} + mk^{n+1})$. 
\end{restatable}

\subsubsection{Width Parameters and Graph Powers} 

As our last main contribution, in \Cref{s-powers}, we address the natural question of whether boundedness of some width parameter is preserved under graph powers. Clearly, this cannot hold for treewidth. However, Suchan and Todinca~\cite{ST07} showed that, for every positive integer $r$ and graph $G$, $\nlcw(G^r) \leq 2(r + 1)^{\nlcw(G)}$, where NLC-width ($\nlcw$) is a parameter equivalent to clique-width. Jaffke et al.~\cite{JKST19} showed that, for every positive integer $r$ and graph $G$, $\mimw(G^r) \leq 2\mimw(G)$. Lima et al.~\cite{LMM24} showed that, for every odd positive integer $r$ and graph $G$, $\tin(G^r) \leq \tin(G)$ and that, for every fixed even positive integer $r$, there is no function $f$ such that $\tin(G^r) \leq f(\tin(G))$ for all graphs $G$. It follows from a result of Bonnet et al.~\cite[Theorem~8.1]{BKTW21} (see also \cite{BFLP24}) that there exists a function $f$ such that, for every positive integer $r$ and graph $G$, $\tww(G^r) \leq f(\tww(G), r)$. Moreover, it is not difficult to see that there cannot be any function uniformly bounding the twin-width of graph powers, i.e. a function $f$ such that, for every positive integer $r$ and graph $G$, $\tww(G^r) \leq f(\tww(G))$. This follows from the fact that the class of leaf powers (induced subgraphs of powers of trees) contains that of interval graphs \cite{BHMW10} and the latter class has unbounded twin-width \cite{BGKTW21}.

We show that sim-width behaves similarly to tree-independence number with respect to graph powers: 

\begin{restatable}{theorem}{powersim}
\label{powersim}
Let $r \geq 1$ be an odd integer and let $G$ be a graph. If $(T,\delta)$ is a branch decomposition of $G$ with $\simw_G(T,\delta) = w$, then $(T,\delta)$ is also a branch decomposition of $G^r$ with $\simw_{G^r}(T,\delta) \leq w$. In particular, $\simw(G^r) \leq \simw(G)$, for every odd integer $r \geq 1$.
\end{restatable}

Moreover, for every fixed even positive integer $r$, we observe that there is no function $f$ such that $\simw(G^r) \leq f(\simw(G))$ for all graphs $G$ (\Cref{p-new}). 
This helps to increase our understanding of the currently poorly understood width parameter sim-width.

\subsection{Consequences of \Cref{t-main,powersim}}
In this section we highlight some algorithmic consequences of \Cref{t-main,powersim}. Before doing so, we first need to briefly review the algorithmic applications of tree-independence number and the more powerful sim-width. 

Concerning tree-independence number, Lima et al.~\cite{LMM24} showed that, for each even positive integer $d$, \textsc{Max Weight Distance-$d$ Packing} can be solved in $\mathsf{XP}$ time parameterized by the independence number of a given tree decomposition of the input graph. This problem is a common generalization of \textsc{Max Weight Distance-$d$ Independent Set} and \textsc{Max Weight Induced Matching}, among others. They also obtained an algorithmic meta-theorem for the problem of finding a maximum-weight induced subgraph with bounded chromatic number satisfying an arbitrary but fixed property expressible in counting monadic second-order logic ($\mathsf{CMSO}_2$) which generalizes \textsc{Feedback Vertex Set} and \textsc{Odd Cycle Transversal}, among others. Dallard, Milani\v{c} and \v{S}torgel~\cite{DMS24} showed that \textsc{$k$-Clique} and \textsc{List $k$-Colouring} admit linear-time algorithms and that \textsc{List $k$-Edge Colouring} admits a quadratic-time algorithm, for every graph class of bounded tree-independence number. 

There exists a wealth of algorithmic applications of mim-width (see, e.g., \cite{BV13,BK19,BTV13,GMR20,JKST19,JKT,JKT20}) and a recent meta-theorem provided by Bergougnoux, Dreier and Jaffke~\cite{BDJ22}, which captures all these results and shows that all problems expressible in $\mathsf{A\& C \ DN}$ logic\footnote{An extension of existential $\mathsf{MSO}_1$ logic.} can be solved in $\mathsf{XP}$ time parameterized by the mim-width of a given branch decomposition of the input graph. However, much less is known about algorithmic applications of the more powerful sim-width. To the best of our knowledge, the only $\mathsf{NP}$-hard problem known to be in $\mathsf{XP}$ parameterized by the sim-width of a given branch decomposition of the input graph is {\sc List $k$-Colouring} \cite{MY23}. In \cite{MY23} it is also shown that, if \textsc{Independent Set} is in $\mathsf{XP}$ parameterized by the sim-width of a given branch decomposition of the input graph, then the same is true for its generalization \textsc{Independent $\mathcal{H}$-packing}. The best known result in this direction is a recent quasipolynomial-time algorithm for \textsc{Independent Set} on graphs of bounded sim-width \cite{BKR23}. 

With the aid of \Cref{t-main,powersim}, we can extend the list of $\mathsf{NP}$-hard problems which are in $\mathsf{XP}$ parameterized by the sim-width of a given branch decomposition of the input graph. Consider first \textsc{List $k$-Edge Colouring}, where an instance consists of a graph $G$ and a list of colours $L(e) \subseteq \{1,\ldots, k\}$ for each $e \in E(G)$, and the task is to determine whether there is an assignment of colours to the edges of~$G$ using colours from the lists in such a way that no two adjacent edges receive the same colour. This problem admits a quadratic-time algorithm on every class of bounded tree-independence number, with no tree decomposition required as input \cite{DMS24}. \Cref{t-main} implies the following:

\begin{corollary} For every $k \geq 1$, {\sc List $k$-Edge Colouring} is quadratic-time solvable on every class of bounded sim-width. 
\end{corollary}

Indeed, we simply check whether the input graph of sim-width at most $c$ contains a vertex of degree at least $k+1$. If it does, we have a no-instance. Otherwise, the input graph is $K_{k+1,k+1}$-subgraph-free and, by \Cref{c-main}, has tree-independence number at most $f(c,k)$, for some function $f$ (we remark that all bounding functions from \Cref{c-main} are computable). We then apply the algorithm from \cite{DMS24}.

We now consider \textsc{List $(d,k)$-Colouring}. A \textit{$(d,k)$-colouring} of a graph $G$ is an assignment of colours to the vertices of $G$ using at most $k$ colours such that no two vertices at distance at most $d$ receive the same colour. For fixed $d$ and $k$, \textsc{$(d,k)$-Colouring} is the problem of determining whether a given graph $G$ has a $(d, k)$-colouring. The \textsc{List $(d,k)$-Colouring} problem requires in addition that every vertex $u$ must receive a colour from some given set $L(u) \subseteq \{1,\ldots, k\}$. Clearly, a graph $G$ has a $(d, k)$-colouring if and only if $G^d$ has a $(1,k)$-colouring. Sharp~\cite{Sharp07} provided the following complexity dichotomy: For fixed $d \geq 2$, \textsc{$(d,k)$-Colouring} is polynomial-time solvable for $k \leq \lfloor\frac{3d}{2}\rfloor$ and $\mathsf{NP}$-hard for $k > \lfloor\frac{3d}{2}\rfloor$.

We recall that for every $k \geq 1$, \textsc{List $(1, k)$-Colouring} is in $\mathsf{XP}$ parameterized by the sim-width of a given branch decomposition of the input graph \cite{MY23}. Hence, \Cref{powersim} immediately implies the following:

\begin{corollary}\label{colsim} For every $k \geq 1$ and odd $d \geq 1$, {\sc List $(d,k)$-Colouring} is in $\mathsf{XP}$ parameterized by the sim-width of a given branch decomposition of the input graph.  
\end{corollary}

Kratsch and M{\"{u}}ller~\cite{KM12} showed that \textsc{List $(1, k)$-Colouring} is polynomial-time solvable for AT-free graphs and hence for the subclass of cocomparability graphs. Moreover, Chang et al.~\cite{CHK03} showed that if $G$ is an AT-free graph, then $G^d$ is a cocomparability graph for any $d \geq 2$ (see also \cite{BDXY06}). Therefore, for any $k,d \geq 1$, \textsc{List $(d,k)$-Colouring} is polynomial-time solvable for AT-free graphs. Since we can compute in linear time a branch decomposition of a cocomparability graph with sim-width at most $1$~\cite{KKST17}, \Cref{colsim} implies the following special case:

\begin{corollary} For every $k \geq 1$ and odd $d \geq 1$, {\sc List $(d,k)$-Colouring} is polynomial-time solvable for AT-free graphs. 
\end{corollary}

Finally, we note that the existence of a polynomial-time algorithm for \textsc{Independent Set} for graph classes whose sim-width is bounded and ``quickly'' computable, pipelined with \Cref{powersim} and \cite[Corollary~7]{MY23}, would imply the same for the more general \textsc{Max Weight Distance-$d$ Packing}, for even $d\geq 2$. This would be in contrast to the situation for odd $d \geq 3$, as for any such $d$, \textsc{Max Weight Distance-$d$ Packing} is $\mathsf{NP}$-complete even for chordal graphs~\cite{EGM14}, which have sim-width~$1$~\cite{KKST17}. 

\subsection{Consequences of \Cref{t-main3}}

The {\it subdivision} of an edge $e=uv$ of a graph replaces $e$ by a new vertex $z$ and edges $uz$ and $zv$. 
The {\it claw} is the $4$-vertex star with vertices $u,v_1,v_2,v_3$ and edges $uv_1$, $uv_2$ and $uv_3$.
A {\it subdivided} claw is a graph obtained from a claw by subdividing each of its edges zero or more times. 
The {\it disjoint union} $G_1+G_2$ of two vertex-disjoint graphs $G_1$ and $G_2$ is the graph $(V(G_1)\cup V(G_2),E(G_1)\cup E(G_2))$.
The graph class ${\cal S}$ consists of all non-empty  disjoint unions of a set of zero or more subdivided claws and paths (see Figure~\ref{f-st} for an example of a graph that belongs to ${\cal S}$).
We let ${\cal T}$ consists of all line graphs of graphs in~${\cal S}$.

 \begin{figure}
\begin{minipage}[c]{0.5\textwidth}
\hspace*{1cm}
\begin{tikzpicture}[scale=1]
\draw (0,1)--(-1,1)--(-2,0)--(-1,-1)--(2,-1) (-2,0)--(1,0); \draw[fill=black] (-1,1) circle [radius=2pt] (0,1) circle [radius=2pt] (-2,0) circle [radius=2pt] (-1,0) circle [radius=2pt] (0,0) circle [radius=2pt] (1,0) circle [radius=2pt] (-1,-1) circle [radius=2pt] (0,-1) circle [radius=2pt] (1,-1) circle [radius=2pt] (2,-1) circle [radius=2pt];
\end{tikzpicture}
\end{minipage}
\qquad
\hspace*{-1.5cm}
\begin{minipage}[c]{0.5\textwidth}
\begin{tikzpicture}[scale=1] 
\draw (-1,1)--(0,1) 
(-1,-1)--(2,-1) (-1,0)--(1,0); \draw[fill=black] (-1,1) circle [radius=2pt] (0,1) circle [radius=2pt]  (-1,0) circle [radius=2pt] (0,0) circle [radius=2pt] (1,0) circle [radius=2pt] (-1,-1) circle [radius=2pt] (0,-1) circle [radius=2pt] (1,-1) circle [radius=2pt] (2,-1) circle [radius=2pt];
\end{tikzpicture}
\end{minipage}
\caption{The graph $S_{2,3,4}+P_2+P_3+P_4$, which is an example of a graph that belongs to ${\cal S}$.}\label{f-st}
\end{figure}

Dabrowski, Johnson and Paulusma~\cite{DJP19} showed that for any finite set of graphs ${\cal H} = \{H_1,\ldots,H_k\}$, the class of ${\cal H}$-free line graphs has bounded clique-width if and only if $H_i\in {\cal T}$ for some $i \in \{1,\ldots, k\}$. By using Theorem~\ref{t-main3} we can extend this result as follows.

\begin{corollary}\label{c-line}
For any finite set of graphs ${\cal H} = \{H_1,\ldots, H_k\}$, the following are equivalent:
\begin{itemize}
  \item The class of ${\cal H}$-free line graphs has bounded clique-width.
  \item The class of ${\cal H}$-free line graphs has bounded mim-width.
  \item The class of ${\cal H}$-free line graphs has bounded sim-width.
  \item The class of ${\cal H}$-free line graphs has bounded tree-independence number.
  \item $H_i\in {\cal T}$ for some $i \in \{1,\ldots, k\}$.
\end{itemize}
\end{corollary}

\subsection{Consequences of \cref{tinmim}}\label{conseqtinmim}

In this section we highlight some structural consequences of \Cref{tinmim}. Its algorithmic connections with the problem of computing tree decompositions of small independence number will be deferred to \Cref{s-con}, as mostly contingent upon some open problems.  

Considerable attention has been recently devoted to understanding the substructures of graphs with large treewidth or large pathwidth (a parameter dominated by treewidth). While under the minor and subgraph relations these substructures are well understood thanks to the Grid-minor theorem \cite{RS86} (see \Cref{gridminor}) and the Forest-minor theorem \cite{RS83}, much less is known for the induced subgraph relation.

Most results on induced substructures of graphs with large pathwidth or large treewidth deal with specific graph classes such as classes of bounded degree or defined by finitely many forbidden induced subgraphs (see, e.g., \cite{ACDH24,Hic23,Kor23,LR22} and references therein). For other width parameters, the situation is even more obscure. Given a width parameter $p$, one would like to characterize the families $\mathcal{F}_p$ of unavoidable induced subgraphs of graphs with large $p$. More formally, $\mathcal{F}_p$ is any set of graphs for which there exists an integer $c \in \mathbb{N}$ such that every graph $G$ with $p(G) > c$ contains a member of $\mathcal{F}_p$ as an induced subgraph\footnote{The families ${\cal F}_{\tw}$ are called \textit{useful families} in \cite{ACDH24}.}.

Even though characterizing the families $\mathcal{F}_p$ for fixed $p \in \{\cw, \mimw, \simw, \tin\}$ is widely open, there are some obvious graphs that any $\mathcal{F}_p$ must contain. For example, for each of these four parameters $p$, any $\mathcal{F}_p$ must contain an induced subgraph of every subdivision of a wall and an induced subgraph of the line graph of every subdivision of a wall \cite{BHMPP22,DJP19,DMS24,GW07,KKST17}.
Moreover, any $\mathcal{F}_{\tin}$ must in addition contain a complete bipartite graph~\cite{DMS24}.
Observe also that \Cref{c-line} can be rephrased as follows. For each $p \in \{\cw,\mimw,\simw,\tin\}$ and finitely defined\footnote{A hereditary class is \textit{finitely defined} if the set of its minimal forbidden induced subgraphs is finite.} class $\mathcal{C}$ of line graphs, the unavoidable induced subgraphs of graphs in $\mathcal{C}$ with large $p$ are precisely the graphs in $\mathcal{T}$.  

\Cref{tinmim} readily implies that a graph with large tree-independence number either contains a large complete bipartite graph as an induced subgraph or has large mim-width, and so any $\mathcal{F}_{\tin}$ contains precisely some complete bipartite graph and graphs from some $\mathcal{F}_{\mimw}$. 

\begin{corollary} For every integer $k \geq 1$ and graph $G$ with $\tin(G) \geq 6(2^{2k-1} + k^{k+2})$, either 
\begin{itemize}
\item $G$ contains a $K_{k,k}$ as an induced subgraph, or
\item $\mimw(G) \geq k$.
\end{itemize}
\end{corollary}

\Cref{tinmim} has another structural consequence, related to a conjecture of Dallard, Milani\v{c} and \v{S}torgel~\cite{DMS24III} (see also a note at the end of \Cref{s-con}). A graph class $\mathcal{G}$ is \textit{$(\tw, \omega)$-bounded} if there exists a function $f$ such that, for every $G \in \mathcal{G}$ and induced subgraph $H$ of $G$, $\tw(H) \leq f(\omega(H))$. Ramsey's theorem implies that in every graph class of bounded tree-independence number, the treewidth is bounded by an explicit polynomial function of the clique number, and hence bounded tree-independence number implies $(\tw,\omega)$-boundedness \cite{DMS24}. In fact, a partial converse is conjectured to hold \cite[Conjecture~8.5]{DMS24III}:

\begin{conj}[Dallard, Milani\v{c} and \v{S}torgel~\cite{DMS24III}]\label{twomegaconj} A hereditary graph class is $(\tw, \omega)$-bounded if and only if it has bounded tree-independence number.
\end{conj}

Dallard, Milani\v{c} and \v{S}torgel~\cite{DMS24III} showed that the conjecture holds for every graph class obtained by excluding a single graph $H$ with respect to any of the following containment relations: subgraph, topological minor, minor, and their induced variants. Very recently, Abrishami et al.~\cite{AAC24} showed that it holds for (even hole, diamond, pyramid)-free graphs. We use \Cref{tinmim} and the fact that a $(\tw, \omega)$-bounded graph class is $K_{t,t}$-free for some $t$ \cite{DMS21a} to prove that  \Cref{twomegaconj} holds for any (not necessarily hereditary) graph class of bounded mim-width.

\begin{corollary}\label{cor-conj} A graph class of bounded mim-width is $(\tw, \omega)$-bounded if and only if it has bounded tree-independence number.
\end{corollary}

Note that there exist $(\tw,\omega)$-bounded graph classes of unbounded mim-width, for example chordal graphs or even the proper subclass of split graphs~\cite{KKST17}.

\section{Preliminaries}\label{s-pre}

\textcolor{black}{We only consider graphs that are {\it simple} (i.e., with no self-loops and no multiple edges), unless otherwise stated.} Since we will not directly work with clique-width and twin-width, we refer the reader to \cite{DJP19} and \cite{BKTW21}, respectively, for their definitions. In some cases, it will be convenient to work with parameters equivalent to treewidth on the class of all graphs (branch-width and mm-width) and to recall that clique-width is equivalent to rank-width. For these reasons, we provide the definitions of branch-width, mm-width and rank-width as well. For all standard graph-theoretic notions not defined here, we refer the reader to \cite{Die}. 

Given a set $S$, a function $f\colon 2^S \rightarrow \mathbb{Z}$ is \emph{symmetric} if $f(X) = f(\overline{X})$ for all $X \subseteq S$, where we use $\overline{X}$ to denote $S\setminus X$.
A \emph{branch decomposition on $S$} is a pair $(T, \delta)$, where $T$ is a subcubic tree and $\delta$ is a bijection between $S$ and the leaves of $T$. Each edge $e \in E(T)$ naturally splits the leaves of the tree in two parts depending on what component they belong to when $e$ is removed. In this way, each edge $e \in E(T)$ represents a partition of $S$ into two partition classes that we denote by $A_{e}$ and $\overline{A_{e}}$.
Let $f\colon 2^S \rightarrow \mathbb{Z}$ be a symmetric function and let $(T, \delta)$ be a branch decomposition on $S$.
The \textit{$f$-width of $(T, \delta)$} is the quantity $\max_{e \in E(T)}f(A_e)$.
The \textit{$f$-branch-width on $S$} is either the minimum $f$-width over all branch decompositions on $S$
when $|S| \ge 2$, or $f(\varnothing)$ when $|S| \le 1$.

Let $G = (V, E)$ be a graph.
For $X \subseteq E$, let $\mathrm{mid}(X)$ be the set of vertices that are incident with both an edge in $X$ and another edge in $E\setminus X$, and let $\eta_G(X) = |\mathrm{mid}(X)|$.
We define the \textit{branch-width of $G$}, denoted $\bw(G)$, to be the $\eta_G$-branch-width on $E$.
For $X \subseteq V$, let $G[X, \overline{X}]$ denote the edges of $G$ having one endpoint in $X$ and the other endpoint in $\overline{X}$, and let $\mathrm{cutmm}_{G}(X)$ be the size of a maximum matching in $G[X, \overline{X}]$.
The \textit{mm-width of $G$}, denoted $\mmw(G)$, is the $\mathrm{cutmm}_{G}$-branch-width on $V$.
For $X \subseteq V$, let $\mathrm{cutmim}_{G}(X)$ be the size of a maximum induced matching in $G[X, \overline{X}]$.
The \textit{mim-width of $G$}, denoted $\mimw(G)$, is the $\mathrm{cutmim}_{G}$-branch-width on $V$.
For $X \subseteq V$, let $\mathrm{cutsim}_{G}(X)$ be the size of a maximum induced matching between $X$ and $\overline{X}$ in $G$.
The \textit{sim-width of $G$}, denoted $\simw(G)$, is the $\mathrm{cutsim}_{G}$-branch-width on $V$.
For a bipartite graph $G$ with bipartition $X$ and $\overline{X}$, its bipartite adjacency matrix is an $|X| \times |\overline{X}|$ matrix $\mathbf{B}_{G} = (b_{i,j})_{i\in X, j\in \overline{X}}$ over the binary field $\mathbf{GF}(2)$ such that $b_{i,j} = 1$ if and only if $\{i, j\} \in E(G)$. For $X \subseteq V$, let $\mathrm{cutrk}_{G}(X) = \mathrm{rank} \ \mathbf{B}_{G[X, \overline{X}]}$.
The \textit{rank-width of $G$}, denoted $\rw(G)$, is the $\mathrm{cutrk}_{G}$-branch-width on $V$.
When it is clear from context, we refer to a branch decomposition on $V(G)$ or $E(G)$ as a branch decomposition of $G$.

A \textit{tree decomposition} of a graph $G$ is a pair $\mathcal{T} = (T, \{X_t\}_{t\in V(T)})$, where $T$ is a tree such that every node $t$ of $T$ is assigned a vertex subset $X_t \subseteq V(G)$, called a \textit{bag}, satisfying the following three conditions: 
\begin{description}
\item[(T1)] every vertex of $G$ belongs to at least one bag; 
\item[(T2)] for every edge $uv \in E(G)$, there exists a bag containing both $u$ and $v$;
\item[(T3)] for every vertex $u \in V(G)$, the subgraph $T_u$ of $T$ induced by $\{t \in V(T) : u \in X_t\}$ is connected. 
\end{description}
\noindent
The \textit{width} of $\mathcal{T}$ is the maximum value of $|X_t| - 1$ over all $t \in V(T)$. The \textit{treewidth} of a graph $G$, denoted $\tw(G)$, is the minimum width over all tree decompositions of $G$. Given a graph $G$ and a tree decomposition $\mathcal{T} = (T, \{X_t\}_{t\in V(T)})$ of $G$, the \textit{independence number of $\mathcal{T}$}, denoted $\alpha(\mathcal{T}$), is the quantity $\max_{t\in V(T)} \alpha(G[X_t])$, where $\alpha(H)$ denotes the maximum size of an independent set in a graph $H$. The \textit{tree-independence number} of $G$, denoted $\tin(G)$, is the minimum independence number over all tree decompositions of $G$.

We already explained in \Cref{sec:intro} the relationships between the six parameters illustrated in \Cref{f-width}. In the following, we recall the results showing that a certain parameter $p$ dominates another parameter $q$ and the corresponding computable functions. Robertson and Seymour~\cite{RS91} showed that, for every graph $G$, 
\begin{equation}\label{twequiv}
  \bw(G)-1 \le \tw(G) \le \left\lfloor\frac{3}{2}\bw(G)\right\rfloor-1,
\end{equation}
so treewidth and branch-width are linearly equivalent. More recently, Vatshelle~\cite{Vat12} and Jeong et al.~\cite{JST18} showed that, for every graph~$G$,
\begin{equation}
  \label{mmwequiv}
  \mmw(G) \le \bw(G) \le \tw(G)+1 \le 3\mmw(G),
\end{equation}
so mm-width is also linearly equivalent to treewidth and branch-width. Oum and Seymour~\cite{OS06} showed that, for every graph~$G$, 
\begin{equation}
  \label{cwequiv}
  \rw(G) \le \cw(G) \le 2^{\rw(G)+1}-1,
\end{equation}
so clique-width and rank-width are equivalent. Note also that, for every graph~$G$,
\begin{equation}
  \label{widthcontainments1}
  \simw(G) \leq \mimw(G) \le \rw(G) \le \bw(G),
\end{equation}
where the first inequality follows from the definition, the second follows from \cite[Lemma~2.4]{BK19}, for example, and the third is due to Oum~\cite{Oum08}. \Cref{twequiv,cwequiv,widthcontainments1} imply that clique-width dominates treewidth. This was also shown by Corneil and Rotics~\cite{CR05}: more precisely, for every graph~$G$,
\begin{equation}
  \label{widthcontainments2}
  \cw(G) \le 3 \cdot 2^{\tw(G)-1}.
\end{equation}
Munaro and Yang~\cite{MY23} showed that sim-width dominates tree-independence number, and Yolov~\cite{Yol18} and Dallard et al.~\cite{DMS24} showed that tree-independence number in turn dominates treewidth: more precisely, for every graph~$G$,
\begin{equation}
  \label{tinsim}
  \simw(G) \leq \tin(G) \le \tw(G) + 1.
\end{equation}
The following result of Bonnet et al.~\cite{BKTW21} shows that twin-width dominates clique-width: for every graph $G$,
\begin{equation}\label{twinwidthcontainments}
  \tww(G) \le 2^{\cw(G)+1}-1.
\end{equation}

We conclude this section with two results which will be repeatedly used in our proofs but first require a series of definitions. Let $K_t \boxminus K_t$ be the graph obtained from $2K_t$ by adding a perfect matching between the two copies of $K_t$, and let $K_t \boxminus S_t$ be the graph obtained from $K_t\boxminus K_t$ by removing all the edges in one of the complete graphs. We write $R(s,t)$ to denote the minimum number such that any graph on at least $R(s,t)$ vertices contains either a clique of size $s$, or an independent set of size $t$ (such a number exists due to Ramsey's theorem). A \emph{wall of height $h$ and width $r$} (an \textit{$h \times r$-wall} for short) is the graph that can be obtained from the grid of height~$h$ and width~$2r$ as follows. Let $C_1, \dots, C_{2r}$ be the sets of vertices in each of the $2r$ columns of the grid, ordered from left to right. For each column $C_{j}$, let $e_{1}^{j}, e_{2}^{j}, \dots, e_{h-1}^{j}$ be the edges between two vertices of $C_j$, in top-to-bottom order. If $j$ is odd, delete all edges $e_{i}^{j}$ with $i$ even. If $j$ is even, delete all edges $e_{i}^{j}$ with $i$ odd. By removing all degree-$1$ vertices of the resulting graph, we obtain an \textit{elementary $h \times r$-wall}. Finally, an $h \times r$-wall is any subdivision of the elementary $h \times r$-wall.  See Figure~\ref{fig:walls1} for a small example. A \emph{net-wall} is a graph obtained from a wall~$G$ by replacing every vertex~$u$ of $G$ that has three distinct neighbours $v,w,x$ by three new vertices $u_v$, $u_w$, $u_x$ and edges $u_vv$, $u_ww$, $u_xx$, $u_vu_w$, $u_vu_x$, $u_wu_x$. See again Figure~\ref{fig:walls1}.

\begin{figure}[tb]
\begin{minipage}{0.5\textwidth}
\centering
\begin{tikzpicture}[scale=0.85]
\draw (-2,0)--(2,0)--(2,1)--(-3,1)--(-3,2)--(1,2)--(1,1) (-2,0)--(-2,1) (0,0)--(0,1) (-1,1)--(-1,2);
\draw[fill=black]
(-2,0) circle [radius=2pt] (-1,0) circle [radius=2pt] (0,0) circle [radius=2pt] (1,0) circle [radius=2pt] (2,0) circle [radius=2pt] (-3,1) circle [radius=2pt] (-2,1) circle [radius=2pt] (-1,1) circle [radius=2pt] 
(0,1) circle [radius=2pt] (1,1) circle [radius=2pt] (2,1) circle [radius=2pt] (-3,2) circle [radius=2pt] (-2,2) circle [radius=2pt] (-1,2) circle [radius=2pt] (0,2) circle [radius=2pt] (1,2) circle [radius=2pt];
\end{tikzpicture}
\end{minipage}
\begin{minipage}{0.2\textwidth}
\centering
\begin{tikzpicture}[scale=0.85]
\draw(2,0)--(-2,0)--(-2,0.7)--(-2.3,1)--(-3,1)--(-3,2)--(1,2)--(1,1.3)--(1.3,1)--(2,1)--(2,0) (-1,1.7)--(-1.3,2) (-0.7,2)--(-1,1.7)--(-1,1.3)--(-1.3,1)--(-1.7,1)--(-2.3,1) (-1.7,1)--(-2,0.7)
(-0.3,0)--(0,0.3)--(0,0.7)--(-0.3,1)--(-0.7,1)--(-1.3,1) (-0.7,1)--(-1,1.3) (-0.3,1)--(0.3,1)--(0.7,1)--(1,1.3) (0.7,1)--(1.3,1) (0.3,1)--(0,0.7) (0,0.3)--(0.3,0);
\draw[fill=black] (-2,0) circle [radius=2pt] (-1,0) circle [radius=2pt] (-0.3,0) circle [radius=2pt] (0.3,0) circle [radius=2pt] (0,0.3) circle [radius=2pt] (1,0) circle [radius=2pt] (2,0) circle [radius=2pt] 
(-3,1) circle [radius=2pt] (-1.7,1) circle [radius=2pt] (-2.3,1) circle [radius=2pt] (-2,0.7) circle [radius=2pt] (-1.3,1) circle [radius=2pt] (-0.7,1) circle [radius=2pt] (-1,1.3) circle [radius=2pt] (-0.3,1) circle [radius=2pt]
(0.3,1) circle [radius=2pt] (0,0.7) circle [radius=2pt] (0.3,1) circle [radius=2pt] (0.7,1) circle [radius=2pt] (1.3,1) circle [radius=2pt] (1,1.3) circle [radius=2pt] (2,1) circle [radius=2pt] 
(-3,2) circle [radius=2pt] (-2,2) circle [radius=2pt] (-0.7,2) circle [radius=2pt] (-1.3,2) circle [radius=2pt] (-1,1.7) circle [radius=2pt] (0,2) circle [radius=2pt] (1,2) circle [radius=2pt];
\end{tikzpicture}
\end{minipage}
\caption{An elementary wall of height~$2$ and a net-wall.}\label{fig:walls1}
\end{figure}

\begin{proposition}[Kang et al.~\cite{KKST17}]\label{cliquefree} Let $G$ be a graph with $\simw(G) = w$ and no induced subgraph isomorphic to $K_t \boxminus K_t$ and $K_t \boxminus S_t$. Then $\mimw(G) \leq R(R(w + 1, t), R(t, t))$. 
\end{proposition}

\begin{corollary}\label{simwall} There exists a function $f\colon \mathbb{N} \rightarrow \mathbb{N}$ such that the sim-width of the elementary $f(m) \times f(m)$-wall is at least $m$. 
\end{corollary}

\begin{proof} It follows from \Cref{cliquefree} and the fact that the elementary $m\times m$-wall, with $m \geq 7$, has mim-width at least $\sqrt{m}/50$ \cite{BHMPP22}. 
\end{proof}

\section{$\mathbf{K_{t,t}}$-Subgraph-Free Graphs}\label{s-main}

\subsection{The Proof of Theorem~\ref{t-main}}

We begin by considering the class of $K_{t,t}$-subgraph-free graphs. The following result is useful (notice that it follows from a more general result of Gravier et al.~\cite[Theorem~2]{GMRT04}).

\begin{lemma}[{Dabrowski, Demange and Lozin~\cite[Lemma~1]{DDL13}}]
  \label{nokttmatching}
  For any positive integers $t$ and $p$, there exists a number $N(t,p)$ such that every bipartite graph with a matching of size at least $N(t,p)$ and having no $K_{t,t}$-subgraph contains an induced matching of size $p$.
\end{lemma}

\begin{proposition}\label{Kttsubgraph} Let $t$ be a positive integer. Then sim-width, mim-width, tree-independence number, clique-width and treewidth are equivalent for $K_{t, t}$-subgraph-free graphs.
\end{proposition}

\begin{proof}
  Define $f_t\colon \mathbb{N} \rightarrow \mathbb{N}$ such that $f_t(x)$ is the least integer such that every bipartite graph with a matching of size at least $f_t(x)$ and having no $K_{t,t}$-subgraph contains an induced matching of size $x$ (such a number exists by \cref{nokttmatching}).
  Observe that $f_t(x)$ is non-decreasing, and define $f_t^{-1}\colon \mathbb{N} \rightarrow \mathbb{N}$ such that $f_t^{-1}(y) = x$ when $f_t(x) \le y$ but $f_t(x+1) > y$. Note that $f_{t}^{-1}$ is non-decreasing as well.

  Suppose $\mmw(G) = k$.
  Let $p = f_t^{-1}(k)$, so $p$ is an integer such that $k \ge f_t(p)$.
  Then, for any branch decomposition $(T,\delta)$ of $G$, there is some $e \in E(T)$ such that $\mmw_G(A_e) \geq k$.
  Since $G[A_e,\overline{A_e}]$ is a bipartite graph with no $K_{t,t}$-subgraph, the definition of $f_t$ implies that $\cutmim_G(A_e) \ge p$.
  Thus $\mimw(G) \ge p = f^{-1}_t(\mmw(G))$.
  As $\mimw(G) \le \mmw(G)$ by definition, mmw-width and mim-width (and hence also treewidth and clique-width) are equivalent for the class of graphs with no $K_{t,t}$-subgraph. 
  
Observe now that if a graph $G$ has no $K_{t,t}$-subgraph, then it has no $K_{2t}$-subgraph and so no induced $K_{2t} \boxminus K_{2t}$ and $K_{2t} \boxminus S_{2t}$ either. Therefore, \Cref{cliquefree} and \Cref{widthcontainments1} imply that sim-width is equivalent to mim-width for the class of $K_{t, t}$-subgraph-free graphs.   

Finally, since tree-independence number always dominates treewidth and sim-width always dominates tree-independence number (see \Cref{f-width}), we conclude that tree-independence number is equivalent to the other parameters for the class of $K_{t, t}$-subgraph-free graphs.
\end{proof}

An analogue to \cref{Kttsubgraph} for $(K_s,K_{t,t})$-free graphs can be proved, mutatis mutandis, using the following result.

\begin{lemma}[{Dabrowski, Demange and Lozin~\cite[Lemma~2]{DDL13}}]
  \label{nokttmatching2}
  For any positive integers $s$, $t$ and $p$, there exists a number $N'(s,t,p)$ such that every $(K_s,K_{t,t})$-free graph with a matching of size at least $N'(s,t,p)$ contains an induced matching of size $p$.
\end{lemma}

\begin{corollary}
  \label{simequivnoksktt}
  Let $s$ and $t$ be positive integers.
  Then tree-independence number, sim-width, mim-width, clique-width and treewidth are equivalent for the class of $(K_s,K_{t, t})$-free graphs.
\end{corollary}

We provide an alternative proof of \Cref{simequivnoksktt} which makes use of \Cref{tinmim} (recall that the latter result will be proved in \Cref{s-main2}).

\begin{proof}[Alternative proof of \Cref{simequivnoksktt}] Let $G$ be a $(K_s, K_{t, t})$-free graph with $\simw(G) = w$. Since $G$ contains no induced subgraph isomorphic to $K_{s} \boxminus K_{s}$ and $K_{s} \boxminus S_{s}$, \Cref{cliquefree} implies that $\mimw(G) \leq R(R(w + 1, s), R(s, s))$. Let $k = R(R(w + 1, s), R(s, s))$. Since $G$ is 
$K_{t,t}$-free, \Cref{tinmim} implies that $\tin(G) \leq 6(2^{t+k-1} + tk^{t+1})$. Since $G$ is $K_{s}$-free, \cite[Lemma~3.2]{DMS24} implies that $\tw(G) \leq R(s,6(2^{t+k-1} + tk^{t+1})+1)-2$. It remains to finally recall that clique-width always dominates treewidth and that sim-width always dominates clique-width (see \Cref{f-width}).
\end{proof}

We are finally ready to prove \Cref{t-main}, which we restate for convenience.

\tmain*

\begin{proof} The equivalence between sim-width, mim-width, clique-width, treewidth and tree-independence number follows from \Cref{simequivnoksktt}. Observe finally that twin-width is not dominated by any of these parameters, even for $(K_3,K_{2,2})$-free graphs. Indeed, walls are $(K_3,K_{2,2})$-free and have bounded twin-width, but each of the other parameters is unbounded, by \cref{simwall}.
\end{proof}

\subsection{Induced Matchings in $\mathbf{d}$-Degenerate Graphs}

Vatshelle~\cite[Lemma~4.3.9]{Vat12} had a purported lemma used to obtain a lower bound on the mim-width of grids, stating that a $d$-degenerate bipartite graph with a matching of size $\mu$ has an induced matching of size at least $\mu / (d+1)$. Recall that a graph is \textit{$d$-degenerate} if each of its subgraphs has a vertex of degree at most $d$. Vatshelle's statement was later used in \cite[Lemma~5.2]{KKST17} and \cite[Lemma~7]{Men18} to give a lower bound for mim-width that was linear in treewidth for $d$-degenerate graphs. We now give a counterexample to Vatshelle's statement and then provide an alternative linear lower bound for mim-width in terms of treewidth for degenerate graphs (in fact, more generally, for graphs of bounded maximum average degree).

Our counterexample consists of a $d$-degenerate bipartite graph with a matching of size $2d$ but whose largest induced matching is of size $1$. This demonstrates that, for a $d$-degenerate bipartite graph with a matching of size $\mu$, at best one could hope to guarantee an induced matching of size $\mu / 2d$.

\begin{example}\label{newcounter} Let $d$ be a positive integer, and let $V_1, V_2, V_3, V_4$ be pairwise disjoint sets with $|V_i| = d$ for each $i \in \{1,2,3,4\}$. Let $G$ be the graph with vertex set
 $V = V_1\cup V_2\cup V_3 \cup V_4$
 and edge set $E = \bigcup_{i \in \{1,2,3\}} \{uv : u \in V_i, v \in V_{i+1}\}$. Clearly, $G$ is bipartite with vertex bipartition $(V_1 \cup V_3, V_2 \cup V_4)$. Consider now an ordering of $V$ in which $V_1$ precedes $V_2$, $V_2$ precedes $V_3$, and $V_3$ precedes $V_4$. Each vertex of $G$ has at most $d$ earlier neighbours in such an ordering. Hence, $G$ is $d$-degenerate. Moreover, $G$ has a matching of size $2d$ but no induced matching of size $2$.    
\end{example}

In order to linearly lower bound mim-width in terms of treewidth for $d$-degenerate graphs, we make use of the following result of Kanj et al.~\cite[Lemma~4.10]{KPSX11}. Recall that the \textit{maximum average degree} of a graph $G$ is the quantity $\max_{H\subseteq G} 2|E(H)|/|V(H)|$, where the maximum is taken over all subgraphs of $G$.

\begin{lemma}[Kanj et al.~\cite{KPSX11}]\label{avgdegreeim} Let $G$ be a graph with maximum average degree at most $d$. If $G$ has a matching of size $\mu$, then $G$ has an induced matching of size at least $\frac{\mu}{2d-1}$.
\end{lemma}

Since a $d$-degenerate graph has maximum average degree at most $2d$ (see, e.g., \cite{NO12}), we have the following consequence of \cref{avgdegreeim}.

\begin{lemma}\label{inducedmatching} Let $G$ be a $d$-degenerate graph with a matching of size $\mu$. Then $G$ has an induced matching of size at least $\frac{\mu}{4d-1}$.
\end{lemma}

The next result is a straightforward consequence of \cref{avgdegreeim} and the equivalence of treewidth and mm-width (see \cref{mmwequiv}). It implies that mim-width and treewidth are equivalent for a class of graphs with bounded maximum average degree (and, in particular, for degenerate graphs, or graphs with bounded maximum degree). 

\begin{lemma}
  For a graph $G$ with maximum average degree at most $d$,
  \[\mimw(G) \geq \frac{\mmw(G)}{2d-1} \geq \frac{\tw(G)+1}{3(2d-1)}.\]
\end{lemma}

\section{Line Graphs}\label{s-main3}

\subsection{The Proof of Theorem~\ref{t-main3}}

Recall that Gurski and Wanke \cite{GW07} showed that for a class of line graphs $\{L(G) : G \in \mathcal{G}\}$, clique-width is equivalent to treewidth for the underlying graph class $\mathcal{G}$. In this section, we show that, in fact, clique-width (and hence treewidth for the underlying graph class) is also equivalent to mim-width, sim-width and tree-independence number.

Due to known results, it suffices to prove that there is a non-decreasing unbounded function $f$ such that $\simw(L(G)) \ge f(\tw(G))$ for every graph $G$; we prove this as \cref{simline}.
Towards this, we require a preliminary result.
In \cite[Lemma~4.5]{KKST17}, it is shown that the sim-width of a graph cannot increase when contracting an edge. In \cref{edgecontraction}, we show that the sim-width of the corresponding line graph cannot increase either.
We first require a definition.

Let $T$ be a tree with at least one vertex of degree at least three, and let $v$ be a leaf of $T$. Let $u$ be the vertex of $T$ with degree at least $3$ and having shortest distance in $T$ from $v$.  Let $P$ be the $v,u$-path in $T$. The operation of \textit{trimming $v$} consists of deleting the vertices $V(P) \setminus \{u\}$ from $T$. Observe that the operation of trimming decreases the number of leaves of $T$ by exactly $1$.   

\begin{lemma}\label{edgecontraction}
Let $G$ be a graph with $|E(G)| \geq 3$ and let $G'$ be the graph obtained by contracting an edge of $G$. Then $\simw(L(G)) \geq\simw(L(G'))$.
\end{lemma}

\begin{proof} Suppose that $G'$ is obtained from $G$ by contracting\footnote{\textcolor{black}{Recall that the operation of contracting an edge $uv$ in $G$ deletes $u$ and $v$ from $G$ and adds a new vertex~$w$ that is made adjacent to the vertices in $(N_G(u)\setminus \{v\}) \cup (N_G(v)\setminus \{u\})$; note that the resulting graph is simple again.}} the edge $uv$ into the vertex $w$. Let $X = \{x_1,\ldots,x_n\}$ be the set of vertices in $G$ adjacent to $u$ but not $v$. Let $Y = \{y_1,\ldots,y_m\}$ be the set of vertices in $G$ adjacent to $v$ but not $u$. Let $Z = \{z_1,\ldots,z_{\ell}\}$ be the set of vertices in $G$ adjacent to both $u$ and $v$. Suppose that $\simw(L(G)) = k$. Then, there exists a branch decomposition $(T,\delta)$ of $L(G)$ such that $\simw_{L(G)}(T,\delta)=k$. We show how to construct a branch decomposition $(T',\delta')$ of $L(G')$ such that $\simw_{L(G')}(T',\delta') \leq k$. This would imply that $\simw(L(G')) \leq \simw_{L(G')}(T',\delta') \leq k = \simw(L(G))$, thus concluding the proof. We may assume that $L(G') \not\cong K_1$, for otherwise $\simw(L(G')) = 0$ and $\simw(L(G)) \ge 0$.

We construct $(T',\delta')$ from $(T,\delta)$ as follows. First, we trim the leaf $\delta(uv)$ of $T$. Then, for each $z_i$, we recursively trim the leaf $\delta(z_iv)$. We call the resulting tree $T'$. We now argue that these operations can indeed be performed. Since we aim to trim $\ell + 1$ leaves, it is enough to show that $T$ contains at least $\ell + 1$ vertices of degree at least $3$. We first recursively contract edges of $T$ having at least one endpoint of degree $2$. The resulting tree~$\widetilde{T}$ has each of its nodes of degree either $1$ or $3$. Moreover, the number of degree-$3$ vertices in $\widetilde{T}$ is the same as that of $T$ and the same holds for degree-$1$ vertices. Let $d_1$ and $d_3$ be the number of degree-$1$ and degree-$3$ vertices in $\widetilde{T}$ (and hence in $T$). It is easy to see that $d_1 - d_3 = 2$. Since the leaves of $T$ are in bijection with the edges of $G$, we obtain that $d_3 = |E(G)| - 2$. Thus, it suffices to show that $|E(G)| \geq \ell + 3$. We now count the edges of $G$. Since $uv \in E(G)$ and, for each $i \in \{1,\ldots, \ell\}$, we have $\{z_iu, z_iv\} \subseteq E(G)$, we see that $|E(G)| = 1 + 2\ell + k$, for some $k \geq 0$.  In particular, as $|E(G)| \geq 3$, either $|E(G)| \geq \ell + 3$, or $k = 0$ and $\ell = 1$. But in the latter case $L(G') \cong K_1$. Therefore $T'$ is a well-defined subcubic tree. To conclude the construction of $(T', \delta')$, we define $\delta'$ as follows. For each $x_i$, let $\delta'(x_iw) = \delta(x_iu)$; for each $y_i$, let $\delta'(y_iw) = \delta(y_iv)$; for each $z_i$, let $\delta'(z_iw) = \delta(z_iu)$; finally, let $\delta'$ coincide with $\delta$ on the remaining vertices of $L(G')$ (these correspond to the edges of $G$ not adjacent to $uv$). 

We now show that $\simw_{L(G')}(T',\delta') \leq k$. Suppose, to the contrary, that $\simw_{L(G')}(T',\delta') \geq k+1$. Then, there exists $e \in E(T')$ such that $\cutsim_{L(G')}(A'_e,\overline{A'_e}) \geq k+1$. Since $T'$ is obtained from $T$ by repeated applications of trimming, and hence by repeated vertex deletions, $e \in E(T') \cap E(T)$. Then, $e$ naturally induces a partition $(A_e,\overline{A_e})$ of $V(L(G))$. Without loss of generality, we assume that $A_e$ agrees with $A'_e$ on the vertices of $L(G)$ corresponding to edges of $G$ not intersecting $\{u,v\}$, and the same for $\overline{A_e}$ and $\overline{A'_e}$. Now, since $\cutsim_{L(G')}(A'_e,\overline{A'_e}) \geq k+1$, there exist independent sets $P'=\{p_1,\ldots,p_{k+1}\} \subseteq A'_e$ and $Q'=\{q_1,\ldots,q_{k+1}\} \subseteq \overline{A'_e}$ of $L(G')$ such that $L(G')[P',Q'] \cong (k+1)P_2$. If none of the vertices in $P' \cup Q'$ correspond to an edge of $G'$ incident to $w$, then $P' \subseteq A_e$ and $Q' \subseteq \overline{A_e}$, and $P'$ and $Q'$ are independent sets of $L(G)$ such that $L(G)[P',Q'] \cong (k+1)P_2$. This implies that $\simw_{L(G)}(T,\delta) \geq k+1$, a contradiction. Hence, there is a vertex in $P' \cup Q'$, say without loss of generality $p_1$, which corresponds to an edge of $G'$ incident to $w$. Moreover, since $p_1$ is anticomplete to $\{p_j,q_j\}$ in $L(G')$, for each $j \geq 2$, the following claim holds:

\begin{subclaim}\label{adj} For each $j \geq 2$, neither of the edges $p_j$ and $q_j$ of $G'$ is adjacent to $p_1$.
\end{subclaim}

We now introduce the following notation. If $p_1 = x_iw$, we let $a = x_i$ and $b = u$. If $p_1 = y_iw$, we let $a = y_i$ and $b = v$. If $p_1 = z_iw$, we let $a = z_i$ and $b = u$. Note that, in each case, $ab$ is an edge of $G$ and hence a vertex of $L(G)$.

\begin{subclaim}\label{Pindep} Let $P = \{ab,p_2,\ldots,p_{k+1}\}$. Then $P \subseteq A_e$ and $P$ is an independent set of $L(G)$. 
\end{subclaim}

\begin{claimproof}[Proof of \Cref{Pindep}] By \Cref{adj}, for each $j \geq 2$, the edge $p_j$ of $G'$ is incident to neither $w$ nor $a$. Moreover, $b \in V(G) \setminus V(G')$. These two facts imply that $P = \{ab,p_2,\ldots,p_{k+1}\} \subseteq E(G)$ and that $P$ is an independent set of $L(G)$. By construction, $\delta'(p_1) = \delta'(aw) = \delta(ab)$, from which $ab \in A_e$ and so $P \subseteq A_e$. 
\end{claimproof}

To complete the proof of \Cref{edgecontraction}, we distinguish two cases.

\medskip
\noindent \textbf{\textsf{Case 1:}} The edge $q_1$ of $G'$ is not incident to $w$.  

By \Cref{Pindep}, $P = \{ab,p_2,\ldots,p_{k+1}\} \subseteq A_e$ and $P$ is an independent set of $L(G)$. Since $q_1$ is not incident to $w$, it must be that $q_1 \in E(G)$. Similarly, by \Cref{adj}, no $q_j$ with $j \geq 2$ is incident to $w$ and so $Q' = \{q_1, \ldots, q_{k+1}\}\subseteq E(G)$. Therefore, $Q' \subseteq \overline{A_e}$ and $Q'$ is an independent set of $L(G)$. We now show that, in $L(G)$, the vertex $ab$ is adjacent to the vertex $q_1$. As the vertex $p_1$ is adjacent to the vertex $q_1$ in $L(G')$, the edge $q_1$ is adjacent to the edge $p_1=aw$. But by assumption, the edge $q_1$ is not incident to $w$, and so $q_1$ is incident to $a$. Hence, in $G$, the edge $q_1$ is adjacent to the edge $ab$ (note that, since $q_1 \in E(G') \cap E(G)$, no endpoint of $q_1$ belongs to $\{u,v\}$) and so, in $L(G)$, the vertex $q_1$ is adjacent to the vertex $ab$. Observe finally that, since $Q' \subseteq V(L(G'))$, none of the edges of $G'$ in $Q'$ are incident to the vertex $b \in V(G) \setminus V(G')$. Combining this with \Cref{adj} and the fact that $ab$ is adjacent to $q_1$ in $L(G)$, we obtain that $L(G)[P,Q'] \cong (k+1)P_2$. Therefore, $\cutsim_{L(G)}(A_e, \overline{A_e}) \geq k+1$ and so $\simw_{L(G)}(T,\delta) \geq k+1$, a contradiction.

\medskip
\noindent \textbf{\textsf{Case 2:}} The edge $q_1$ of $G'$ is incident to $w$. 

Consider $uv \in V(L(G))$. Either $uv \in A_e$ or $uv \in \overline{A_e}$. 

\noindent \textbf{\textsf{Case 2.1:}} $uv \in  \overline{A_e}$. 

By \Cref{Pindep}, $P = \{ab,p_2,\ldots,p_{k+1}\} \subseteq A_e$ is an independent set of $L(G)$. By assumption and \Cref{adj}, $Q = \{uv,q_2,\ldots,q_{k+1}\}\subseteq \overline{A_e}$. Moreover, again by \Cref{adj}, none of the edges $q_j$ with $j \geq 2$ is adjacent to $uv$ in $G$, or else they would be adjacent to $p_1$ in $G'$. Therefore, $Q$ is an independent set of $L(G)$. It is finally easy to see that, in $L(G)$, $ab$ is adjacent to $uv$, no $q_j$ with $j \geq 2$ is adjacent to $ab$, and no $p_j$ with $j \geq 2$ is adjacent to $uv$. Hence, $L(G)[P,Q] \cong (k+1)P_2$. Therefore, $\cutsim_{L(G)}(A_e, \overline{A_e}) \geq k+1$ and so $\simw_{L(G)}(T,\delta) \geq k+1$, a contradiction.

\noindent \textbf{\textsf{Case 2.2:}} $uv \in  A_e$. 

Let $q_1 = cw$, for some $c \in X \cup Y \cup Z$. By construction, $\delta'(q_1) = \delta'(cw) = \delta(cd)$, for some $d \in \{u,v\}$. Since $q_1 \in \overline{A'_e}$, we have that $cd \in \overline{A_e}$. Let $P = \{uv,p_2,\ldots,p_{k+1}\}$ and $Q = \{cd,q_2,\ldots,q_{k+1}\}$. By symmetry, the same proof as in Case 2.1 applies to show that $P\subseteq A_e$ and $Q\subseteq \overline{A_e}$ are independent sets in $L(G)$ such that $L(G)[P,Q] \cong (k+1)P_2$. Therefore, $\cutsim_{L(G)}(A_e, \overline{A_e}) \geq k+1$ and so $\simw_{L(G)}(T,\delta) \geq k+1$, a contradiction.
\end{proof}

In order to prove \cref{simline}, we also require the following two results, the first of which is an easy observation (see \cite[Lemma~4.5]{KKST17}).

\begin{lemma}\label{vertexdeletion} Let $G$ be a graph and $v \in V(G)$. Then $\simw(G) \geq \simw(G-v)$.
\end{lemma}

\begin{theorem}[Grid-minor theorem \cite{RS86}]\label{gridminor} There exists a non-decreasing unbounded function $g\colon \mathbb{N} \rightarrow \mathbb{N}$ such that, for every $k \in \mathbb{N}$, every graph of treewidth at least $k$ contains the $g(k) \times g(k)$-grid as a minor.   
\end{theorem}

\begin{proposition}\label{simline} There exists a non-decreasing unbounded function $f\colon \mathbb{N} \rightarrow \mathbb{N}$ such that, for every graph $G$, $\simw(L(G)) \geq f(\tw(G))$.
\end{proposition}

\begin{proof} Let $G_n$ denote the $n\times n$-grid. Let $G$ be a graph with $\tw(G) = k$ and let $g$ be the function from \Cref{gridminor}. Then $G$ contains $G_{g(k)}$ as a minor. This implies that there exists a finite sequence of operations $\{f_i\}_{i=1}^{m}$, where each $f_i$ is either an edge contraction or an edge or vertex deletion, such that $f_1(f_2(\cdots f_m(G)\cdots)) \cong G_{g(k)}$. Observe now that, for each operation $f_i$, $\simw(L(G)) \geq \simw(L(f_i(G)))$, by \Cref{edgecontraction} if $f_i$ is an edge contraction, or by \Cref{vertexdeletion} if $f_i$ is an edge or vertex deletion. Therefore, $\simw(L(G)) \geq \simw(L(f_1(f_2(\cdots f_m(G)\cdots)))) = \simw(L(G_{g(k)}))$. But since $L(G_{g(k)})$ is $K_5$-free, \Cref{cliquefree} implies that there exists an increasing function $h\colon\mathbb{N} \rightarrow \mathbb{N}$ such that 
\begin{equation}\label{simlinelower}
h(\simw(L(G_{g(k)}))) \geq \mimw(L(G_{g(k)})).
\end{equation}
Recall now the following well-known fact (see, e.g., \cite{HW18}): Given a minimum width tree decomposition of $L(G)$, replacing each edge with both of its endpoints gives a tree decomposition
of $G$ and so $\tw(L(G)) \geq \frac{1}{2}(\tw(G)+1)-1$. Therefore, since $\tw(G_{g(k)}) = g(k)$, we have that $\tw(L(G_{g(k)})) \geq \frac{1}{2}(g(k)+1)-1$. But the line graph of the $g(k) \times g(k)$-grid is $K_{6,6}$-subgraph-free and so, by the proof of \Cref{Kttsubgraph}, there exists a non-decreasing unbounded function $h'\colon\mathbb{N} \rightarrow \mathbb{N}$ such that $\mimw(L(G_{g(k)})) \geq h'(\mmw(L(G_{g(k)}))) \geq h'(\frac{1}{3}(\tw(L(G_{g(k)})) + 1)) \geq h'(\frac{1}{6}(g(k) + 1))$, where the second inequality follows from \Cref{mmwequiv}. Combining this chain with \Cref{simlinelower}, we obtain that there exists a non-decreasing unbounded function $f\colon\mathbb{N} \rightarrow \mathbb{N}$ such that $\simw(L(G_{g(k)})) \geq f(k)$. This concludes the proof. 
\end{proof}

We are finally ready to prove \Cref{t-main3}, which we restate for convenience.

\tmainthree*
\begin{proof} The implications $2 \Rightarrow 3 \Rightarrow 4$ follow from \Cref{f-width}. The implication $4 \Rightarrow 1$ follows from \Cref{simline}. The implication $1 \Rightarrow 2$ follows from \Cref{cwline}, while $5 \Rightarrow 4$ follows again from \Cref{f-width}. The implication $1 \Rightarrow 5$ follows from the proof of \cite[Lemma~2.4]{BGT98} (see also \cite[Theorem~3.12]{DMS24}). We provide the short argument for completeness. To this end, let $G \in \mathcal{G}$ and let $\mathcal{T} = (T, \{X_t\}_{t \in V(T)})$ be a tree decomposition of width at most $k$, for some $k \in \mathbb{N}$. We build a tree decomposition $\mathcal{T}'$ of $L(G)$ as follows: Replace each bag $X_t$ with the set $B_t$ of edges of $G$ incident with a vertex in $X_t$. It is easy to see that, for each $t \in V(T)$, $\alpha(L(G)[B_t]) \leq |X_t| \leq k + 1$. Hence, $\alpha(\mathcal{T}') \leq k + 1$. Therefore, $1, 2, 3, 4, 5$ are equivalent. 

To complete the proof, it suffices to observe that twin-width is not dominated by any of the other parameters and that treewidth does not dominate any of the other parameters for line graphs (see \Cref{f-line}). As for the former, let $\mathcal{G}$ be the class of $1$-subdivisions of walls. Then $L(\mathcal{G})$ contains the class of net-walls, which has unbounded treewidth (see, e.g., \cite{BHMPP22}). However, $L(\mathcal{G})$ consists of planar graphs and hence has bounded twin-width \cite{BKTW21}. As for the latter, let $\mathcal{G}$ be the class of stars. Then $L(\mathcal{G})$ coincides with the class of complete graphs, which has unbounded treewidth, but for which any other parameter is bounded. 
\end{proof}

\subsection{The Proof of \Cref{mimlinebw}}

In this section we prove \Cref{mimlinebw}, which we restate for convenience.

\mimlinebw*

The upper bound follows from the fact that, for any graph $G$, $\mimw(G) \leq \rw(G)$ and $\rw(L(G))\leq \bw(G)$~\cite{Oum08}. We provide a short direct proof for completeness.

\begin{lemma}\label{uppermim} For any graph $G$, $\mimw(L(G)) \leq \bw(G)$. 
\end{lemma}

\begin{proof} Let $(T, \delta)$ be a branch decomposition on $E(G)$ of minimum $\eta_G$-width $\bw(G)$. Clearly, $(T, \delta)$ is a branch decomposition on $V(L(G))$ as well, so it suffices to show that $\cutmim_{L(G)}(A_e) \leq |\mathrm{mid}(A_e)|$ for any $e \in E(T)$. This follows from the fact that an induced matching of size $k$ in $L(G)[A_e, \overline{A_e}]$ provides $k$ distinct vertices in $G$ that are incident with an edge in $A_e$ and another edge in $\overline{A_e}$.
\end{proof}    

We now turn to the lower bound, which is immediately obtained from the following result. 

\begin{proposition}\label{lowermim} Let $n$ be a positive integer and let $G$ be a graph such that $\bw(G) \geq 25n-1$. Then $\mimw(L(G)) \geq n$.
\end{proposition}

\begin{proof} Let $(T, \delta)$ be an arbitrary branch decomposition on $V(L(G))$.  Then, $(T, \delta)$ is also a branch decomposition on $E(G)$. Since $\bw(G) \geq 25n-1$, there exists $e \in E(T)$ such that $\eta_G(A_e) = |\mathrm{mid}(A_e)| \geq 25n-1$. It is then enough to show that $\cutmim_{L(G)}(A_e,\overline{A_e}) \geq n$. 

Let $M = \mathrm{mid}(A_e)$. Hence, $|M| \geq 25n-1$. For every vertex $m \in M \subseteq V(G)$, there exist edges $x \in A_e$ and $y \in \overline{A_e}$ of $G$ such that both $x$ and $y$ are incident with $m$. We can then define two functions, $l$ and $r$, as follows. The function $l$ assigns to each $m \in M$ a vertex $l(m) \in V(G)$ such that $l(m)m \in A_e$. The function $r$ assigns to each $m \in M$ a vertex $r(m) \in V(G)$ such that $r(m)m \in \overline{A_e}$. Note that $l$ and $r$ are not necessarily injective and that, for each $m \in M$, the vertices $m$, $l(m)$ and $r(m)$ are pairwise distinct. Given $M$, $l$ and $r$ as above, a \textit{perfect triple} $(L,D,R)$ is a triple such that $D \subseteq M$, $L = \{l(d):d \in D\}$, $R = \{r(d):d \in D\}$ and $L$, $D$ and $R$ are pairwise disjoint. The \textit{size} of the perfect triple $(L,D,R)$ is $|D|$. 

Observe that if there exists a perfect triple $(L,D,R)$ of size $n$, then $\cutmim_{L(G)}(A_e,\overline{A_e}) \geq n$. Indeed, suppose that $D = \{d_1,\ldots,d_n\}$. By definition, for each $i \in \{1, \ldots, n\}$, we have that $l(d_i)d_i \in A_e$, $r(d_i)d_i \in \overline{A_e}$ and $l(d_i)d_i$ is adjacent to $r(d_i)d_i$ in $L(G)$. Consider now $X = \{l(d_1)d_1,\ldots,l(d_n)d_n\} \subseteq A_e$ and $Y = \{r(d_1)d_1,\ldots,r(d_n)d_n\} \subseteq \overline{A_e}$. Since $(L,D,R)$ is a perfect triple, $L$, $D$ and $R$ are pairwise disjoint. Therefore, for $i \neq j$, we have that $\{l(d_i), d_i\} \cap \{r(d_j), d_j\} = \varnothing$ and so $L(G)[X,Y] \cong nP_2$. Consequently, $\cutmim_{L(G)}(A_e,\overline{A_e}) \geq n$. 

In view of the paragraph above, it is enough to show that there exists a perfect triple $(L,D,R)$ of size at least $n$. Suppose, to the contrary, that every perfect triple has size less than $n$. Let $(L,D,R)$ be a perfect triple of maximum size $|D| = k < n$. Since $M \neq \varnothing$ and for each $m \in M$, $m$, $l(m)$ and $r(m)$ are pairwise distinct, we have that $k \geq 1$. We now consider the reason why $(L,D,R)$ cannot be extended to a larger perfect triple by adding a vertex $m \in M\setminus D$ to $D$, and possibly $l(m)$ to $L$ and $r(m)$ to $R$, if not already in $L$ and $R$, respectively. We have that either $m \in M\setminus D$ cannot be added (Case 1 below), or $m$ can be added and exactly one of $l(m)$ and $r(m)$ can be added (Cases 2 to 5 below), or $m$ can be added and none of $l(m)$ and $r(m)$ can be added (Case 6 below). There are six possible cases:

\begin{enumerate}
\item $m \in L \cup R$;
\item $m \notin L \cup R$, $l(m) \in D$ and $r(m) \not \in L\cup D$;
\item $m \notin L \cup R$, $l(m) \in R$ and $r(m) \not \in L\cup D$;
\item $m \notin L \cup R$, $l(m) \not \in R\cup D$ and $r(m) \in D$;
\item $m \notin L \cup R$, $l(m) \not \in R\cup D$ and $r(m) \in L$;
\item $m \notin L \cup R$, $l(m) \in R\cup D$ and $r(m) \in L\cup D$.
\end{enumerate}

In the following series of claims, we show that each of these cases holds for a small number of vertices of~$M$. 

\begin{subclaim}\label{case1} Case 1 holds for at most $2k$ vertices of $M$. 
\end{subclaim}

\begin{claimproof}[Proof of \Cref{case1}]
This follows from the fact that $|L| \leq |D|$ and $|R| \leq |D|$, from which $|L\cup R| \leq 2k$. 
\end{claimproof}

\begin{subclaim}\label{case2} Cases 2 and~4 each hold for at most $3k$ vertices of $M$. 
\end{subclaim}

\begin{claimproof}[Proof of \Cref{case2}] By symmetry, it suffices to consider Case 2. Suppose, to the contrary, that there exist $3k+1$ vertices $m \in M\setminus (L\cup D \cup R)$ satisfying Case 2. Since $|D| = k$, the pigeonhole principle implies that there exists $d \in D$ and four vertices $m_1,\ldots,m_4 \in M\setminus (L\cup D \cup R)$ such that $l(m_i)=d$ and $r(m_i) \not \in L\cup D$ for each $i\in\{1,\ldots,4\}$. 
Then, we remove $d$ from $D$ and, in case there is no $d' \in D$ with $d' \neq d$ and $l(d') = l(d)$, we also remove $l(d)$ from $L$ and, in case there is no $d' \in D$ with $d' \neq d$ and $r(d') = r(d)$, we further remove $r(d)$ from $R$. 
We claim that we can add $d$ into $L$, at least two of $m_1, \ldots, m_4$ into $D$, and the corresponding $r(m_i)$'s into $R$ to obtain a perfect triple of size larger than $k$. Since $d \notin L\cup R$ and, for each $i \in \{1,\ldots,4\}$, we have $m_i \notin L\cup R\cup D$ and $r(m_i) \notin L\cup D$, the only possible obstacle is that after adding $m_1$ to $D$, and $r(m_1)$ to $R$, we have that for each $m_j$ with $j\neq 1$, either $m_j = r(m_1)$ (so we cannot add $m_j$ to $D$), or $r(m_j) = m_1$ (so we cannot add $r(m_j)$ to $R$). Observe that there exist distinct indices $p,q \in \{2,3,4\}$ such that $m_p \neq r(m_1)$ and $m_q \neq r(m_1)$, or else two vertices of $m_2, m_3, m_4$ coincide. Without loss of generality, $m_2 \neq r(m_1)$ and $m_3 \neq r(m_1)$. We then assume that $r(m_2) = r(m_3) = m_1$, or else we immediately get a perfect triple of size $k+1$. But in this case we can add $m_2$ and $m_3$ into $D$, $m_1$ into $R$, and $d$ into $L$ to obtain a perfect triple of size $k+1$, a contradiction.
\end{claimproof}

\begin{subclaim}\label{case3} Cases~3 and~5 each hold for at most $6k$ vertices of $M$. 
\end{subclaim}

\begin{claimproof}[Proof of \Cref{case3}]
By symmetry, it suffices to consider Case 3. Suppose, to the contrary, that there exist $6k+1$ vertices $m \in M\setminus (L\cup D \cup R)$ satisfying Case 3. Let $S$ be the set of such vertices. For each $b \in R$, let $w_b = |\{d \in D : r(d)=b\}|$ and $w'_b = |\{s \in S : l(s)=b, r(s) \notin L\cup D\}|$. Note that $\sum_{b \in R}w_b = |D| = k$ and $\sum_{b \in R}w'_b = |S| \geq 6k+1$. Therefore, there exists $c \in R$ such that $w'_c \geq 6w_c+1$. Let $p = w_c$ and let $M' = \{s \in S : l(s) =c, r(s) \notin L \cup D\}$. Hence, $|M'| = w_c' \geq 6p+1$ and take an arbitrary subset $M'' \subseteq M'$ of size $6p+1$.

We now claim that there exists $Q \subseteq M''$ of size at least $p+1$ such that $r(q) \notin Q$ for every $q \in Q$. To see this, for each $m \in M''$, let $\mathrm{deg}(m) = |\{m'' \in M'' : r(m'')=m\}|$. Observe that $\sum_{m \in M''}\mathrm{deg}(m) \leq |M''| = 6p+1$. So the number of vertices $m \in M''$ with $\mathrm{deg}(m) \geq 2$ is at most $3p$ and there are at least $3p+1$ vertices $m \in M''$ with $\mathrm{deg}(m) \leq 1$. Let $M^* = \{m \in M'' : \mathrm{deg}(m) \leq 1\}$. Hence, $|M^*| \geq 3p+1$. Since each $m \in M^*$ satisfies $\mathrm{deg}(m) \leq 1$, there is at most one $m'' \in M''$ such that $r(m'')=m$. We now show how to construct $Q \subseteq M^*$ such that, for every $q \in Q$, $r(q) \notin Q$. Iteratively, for each $m \in M^*$, add $m$ into $Q$ and possibly remove the following vertices from $M^*$ (in case they belong to $M^*$): $m,r(m)$ and the vertex $m'' \in M''$ with $r(m'') = m$ if it exists. At each step, we add one vertex into $Q$ and remove at most three vertices from $M^*$. Therefore, we can repeat the step above $p+1$ times in order to obtain $Q$ of size at least $p+1$. By construction, for each $q \in Q$, we have $r(q)\notin Q$, as desired.

Let $Q$ be a set given by the previous paragraph. We move $c$ from $R$ to $L$, remove from $D$ any $d\in D$ with $r(d)=c$ and add $Q$ into $D$. Moreover, we remove from $L$ any $l \in L$ such that no $d \in D$ satisfies $l(d) =l$. Finally, we add $\{r(q) : q \in Q\}$ into $R$. Observe that $Q \subseteq S$ and $S$ is disjoint from $L \cup D \cup R$. Moreover, each $q \in Q \subseteq M'$ satisfies $l(q) = c$ and $r(q) \notin L\cup D$ and, by the previous paragraph, $r(q) \notin Q$. Therefore, we obtain a perfect triple. In this process we have removed $w_c = p$ vertices from $D$ and added $|Q| \geq p+1 $ vertices into $D$, and so we obtained a perfect triple of larger size, a contradiction. This concludes the proof of \Cref{case3}.
\end{claimproof}

\begin{subclaim}\label{case6} Case 6 holds for at most $4k$ vertices of $M$. 
\end{subclaim}

\begin{claimproof}[Proof of \Cref{case6}]
Suppose, to the contrary, that there exist $4k+1$ vertices $m \in M\setminus (L \cup D \cup R)$ satisfying Case 6. Let $Q$ be the set of such vertices. Hence, $|Q| \geq 4k+1$. Consider now the multigraph $H$ with vertex set $V(H)=V(G)$ and such that, for each $q \in Q$, $l(q)r(q)$ is an edge of $H$. Observe that $l(q)\neq r(q)$ for each $q \in Q$, so $H$ does not have loops. It is well known that every loopless multigraph has a bipartite subgraph with at least half of its edges (see, e.g., \cite[Theorem~1.3.19]{West}). Let $H'$ be a bipartite subgraph of $H$, with bipartition $(V_1,V_2)$, and at least $\lceil|E(H)|/2\rceil \geq \lceil(4k+1)/2\rceil = 2k+1$ edges. Each edge of $H'$ is of the form $e_q = l(q)r(q)$ for some $q \in Q$, and either $l(q) \in V_1$ and $r(q) \in V_2$, or $l(q) \in V_2$ and $r(q) \in V_1$. Without loss of generality, at least $k+1$ edges $e_q$ of $H'$ satisfy the former, and let $E_q$ be the corresponding set. Let $D' = \{q \in Q : e_q \in E_q\}$, $L' = \{l(d) : d\in D'\}$ and $R' = \{r(d) : d \in D'\}$. By construction, $L' \cap R' = \varnothing$. Moreover, for each $q \in D'$, since $q \notin L \cup D \cup R$ and $l(q) \in R \cup D$ and $r(q) \in L\cup D$, we have that $D' \cap (L' \cup R') = \varnothing$. Therefore, $(L',D',R')$ is a perfect triple of size $|D'| \geq k+1$, a contradiction.
\end{claimproof}

By the previous series of claims, there are at most $2k+2\cdot 3k+2\cdot 6k + 4k=24k$ vertices of $M \setminus D$ satisfying at least one of the six cases. Since $|M\setminus D| \geq 24n \geq 24(k + 1)$, there is at least one vertex $m \in M\setminus D$ satisfying none of the six cases. We add $m$ into $D$, and possibly $l(m)$ into $L$ and $r(m)$ into $R$ (if not already present) and obtain a perfect triple of size $k+1$, contradicting the maximality of $(L, D, R)$. Therefore, there exists a perfect triple of size at least $n$, thus concluding the proof of \Cref{lowermim}.
\end{proof} 

Note that the proof of \Cref{lowermim} does not seem to be easily adaptable to bound $\simw(L(G))$, as the functions $l$ and $r$ therein are not necessarily injective.  

We conclude with a comment on the tightness of \Cref{mimlinebw}. The upper bound on $\mimw(L(G))$ is tight in the sense that, for any integer $n \ge 2$, if we let $G =K_{1,n}$, then $\mimw(L(G)) = \bw(G) = 1$.
The problem of determining a tight lower bound on $\mimw(L(G))$ in terms of $\bw(G)$ is left open.

\subsection{Sim-Width of $\mathbf{L(K_{n,m})}$ and $\mathbf{L(K_{n})}$}

We begin by determining the exact value of $\simw(L(K_{n,m}))$. For a positive integer $k$, we use the notation $[k] = \{1, \ldots, k\}$. For positive integers $n$ and $m$, the $n \times m$ \textit{rook graph} $R_{n,m}$ is the graph representing the legal moves of a rook on an $n \times m$ chessboard: the vertex set is $[n] \times [m]$, two vertices being adjacent if and only if they have one of the two coordinates in common. Clearly, $R_{n,m}$ is isomorphic to $K_n \square K_m$, the Cartesian product of $K_n$ and $K_m$. 

\begin{observation}[Folklore]\label{rook} $L(K_{n,m})$ is isomorphic to $R_{n,m}$.
\end{observation}

\begin{proposition}\label{simwbiclique}
Let $n$ and $m$ be integers with $6 < n \leq m$. Then $\simw(L(K_{n,m})) = \lceil{\frac{n}{3}}\rceil$.
\end{proposition}

\begin{proof} In view of \Cref{rook}, we show that $\simw(R_{n,m}) = \lceil{n/3}\rceil$. We first need the following definition. For an integer $l \geq 1$, an \textit{$l$-caterpillar} is a subcubic tree $T$ on $2l$ vertices with $V(T) = \{s_1, \ldots, s_l, t_1, \ldots, t_l\}$, such that $E(T) = \{s_it_i : 1 \leq i \leq l\} \cup \{s_is_{i+1} : 1 \leq i \leq l - 1\}$. Note that we label the leaves of an $l$-caterpillar $t_1, t_2, \ldots, t_l$, in this order.

  Let us begin by showing the upper bound $\simw(R_{n,m}) \leq \lceil{n/3}\rceil$. Let $a = \lceil{n/3}\rceil$ and $b = \lfloor{2n/3}\rfloor$. Observe that $1 \leq a < b < n$. Build an $(a+1)$-caterpillar $C_1$ with leaves $p_1,l_1,l_2,\ldots,l_a$, a $(b-a+1)$-caterpillar $C_2$ with leaves $p_2,l_{a+1},l_2,\ldots,l_b$ and an $(n-b+1)$-caterpillar $C_3$ with leaves $p_3,l_{b+1},l_2,\ldots,l_n$. For each $i \in [n]$, build an $(m+1)$-caterpillar $D_i$ with leaves $h_i,l_{i,1},\ldots, l_{i,m}$ and add the edge $h_il_i$. Finally, add a vertex $p_0$ and add edges $p_0p_1,p_0p_2$ and $p_0p_3$. Let $T$ be the resulting tree and let $\delta$ be the function mapping each $(i,j) \in [n] \times [m]$ to the leaf $l_{i,j}$ of $T$. Clearly, $(T,\delta)$ is a branch decomposition of $R_{n,m}$. We now show that $\simw_{R_{n,m}}(T,\delta) \leq \lceil{n/3}\rceil$.

Let $e \in E(T)$ and let $(A_e,\overline{A_e})$ be the corresponding bipartition of $V(R_{n,m})$. Suppose first that $e = h_il_i$ or $e \in E(D_i)$, for some $i\in [n]$. Then, without loss of generality, the first coordinate of each vertex in $A_e$ equals $i$ and so $A_e$ is a clique in $R_{n,m}$. Therefore, $\cutsim_{R_{n,m}}(A_e,\overline{A_e}) \leq 1$. Suppose instead that $e = p_0p_j$ or $e \in E(C_j)$, for some $j \in \{1,2,3\}$. Then, without loss of generality, each vertex in $A_e$ has first coordinate between $1$ and $a = \lceil{n/3}\rceil$, if $j = 1$, or between $a+1 = \lceil{n/3}\rceil+1$ and $b =  \lfloor{2n/3}\rfloor$, if $j = 2$, or between $b+1 = \lfloor{2n/3}\rfloor+1$ and $n$, if $j = 3$. In any case, it is easy to see that there are at most $\lceil{n/3}\rceil$ distinct choices for the first coordinate, and so $\alpha(R_{n,m}[A_e]) \leq \lceil{n/3}\rceil$, from which $\cutsim_{R_{n,m}}(A_e,\overline{A_e}) \leq \lceil{n/3}\rceil$.

Finally, we show that $\simw(R_{n,m}) \geq \lceil{n/3}\rceil$. Let $D = \{(i,i):1 \leq i \leq n\}$. Suppose, to the contrary, that $\simw(R_{n,m}) < \lceil{n/3}\rceil$ and let $(T, \delta)$ be a branch decomposition with sim-width $w < \lceil{n/3}\rceil$. We first show that there exists $e \in E(T)$ such $|D|/3 \leq |A_e \cap D|,|\overline{A_e}\cap D|$. Indeed, by trimming the set of leaves $\delta(V(R_{n,m}) \setminus D)$ of $T$, we obtain a branch decomposition of $R_{n,m}[D]$ with sim-width at most $w$. We then apply \cite[Lemma~2.3]{KKST17}\footnote{Notice that the statement of the Lemma in \cite{KKST17} contains a typo, as $<$ should be replaced by $\leq$. It should read as follows. Let $G$ be a graph, let $w$ be a positive integer, and let $f\colon 2^{V(G)}\rightarrow\mathbb{N}$ be a symmetric function. If $G$ has $f$-width at most $w$, then $V(G)$ admits a bipartition $(A_1, A_2)$ where $f(A_1) \leq w$ and $|V(G)|/3 \leq |A_1|, |A_2| \leq 2|V(G)|/3$.} to the graph $R_{n,m}[D]$ and the obtained branch decomposition. Fix now $e$ as above and suppose that $A_e \cap D = \{(a_1,a_1),\ldots,(a_r,a_r)\}$ and $\overline{A_e} \cap D = \{(b_1,b_1),\ldots,(b_s,b_s)\}$, for pairwise distinct $a_1, \ldots, a_r, b_1, \ldots, b_s$ in $[n]$ with $r,s \geq \lceil{n/3}\rceil$. For each $i \in [\lceil{n/3}\rceil]$, we proceed as follows. If $(a_i,b_i) \in A_e$, let $x_i = (a_i,b_i)$ and $y_i = (b_i,b_i)$. Else, $(a_i,b_i) \in \overline{A_e}$, and let $x_i = (a_i,a_i)$ and $y_i = (a_i,b_i)$. Let $X = \{x_1,\ldots,x_{\lceil{n/3}\rceil}\}$ and $Y = \{y_1,\ldots,y_{\lceil{n/3}\rceil}\}$. Clearly, $X \subseteq A_e$ and $Y \subseteq \overline{A_e}$. Moreover, $X$ is an independent set, since no two of its vertices share a coordinate, and similarly for $Y$. Finally, each $x_i$ is adjacent to $y_i$ and no other $y_j$ with $j \neq i$. Therefore, $\cutsim_{R_{n,m}}(A_e,\overline{A_e}) \geq \lceil{n/3}\rceil$,  contradicting the fact that $(T, \delta)$ has sim-width $w < \lceil{n/3}\rceil$.
\end{proof}

We conclude with some observations related to $\simw(L(K_{n}))$. Since $L(K_{n})$ contains $L(K_{n/2,n/2})$ as an induced subgraph, \Cref{simwbiclique} implies that $\simw(L(K_{n})) \geq \lceil n/6\rceil$, for $n > 12$. Moreover, since $\bw(K_n) = \left\lceil\frac{2n}{3}\right\rceil$ for $n \geq 3$ \cite{RS91}, \Cref{mimlinebw} implies that $\simw(L(G)) \leq \left\lceil\frac{2n}{3}\right\rceil$. 

\begin{lemma}\label{simwclique}
Let $n > 12$ be an integer. Then $\left\lceil\frac{n}{6}\right\rceil \leq \simw(L(K_{n})) \leq \left\lceil\frac{2n}{3}\right\rceil$.
\end{lemma}

We expect the value of $\simw(L(K_{n}))$ to be close to the lower bound $\left\lceil\frac{n}{6}\right\rceil$ and leave its determination as an open problem.

\section{$\mathbf{K_{t,t}}$-Free Graphs: The Proof of Theorem~\ref{t-main2}}\label{s-main2}

In this section, we consider the class of $K_{t,t}$-free graphs and ask which parameters from \Cref{f-width} become equivalent or comparable when restricted to this class. In fact, we answer this question except for one remaining open case (see \Cref{kttfreeopen}), as shown in \Cref{kttrelfig}.

Let us first analyse the pairs of comparable parameters from \cref{f-width}. 

\begin{lemma}\label{comp} Even when restricted to $K_{2,2}$-free graphs:
  \begin{itemize}
    \item Treewidth does not become equivalent to any of the other parameters in \Cref{f-width};
    \item Clique-width does not become equivalent to twin-width, mim-width or sim-width; 
    \item Mim-width does not become equivalent to sim-width. 
  \end{itemize}
\end{lemma}

\begin{proof} The class of cliques is $K_{t,t}$-free for any $t$, has unbounded treewidth but bounded clique-width, twin-width and tree-independence number. Therefore, treewidth does not become equivalent to any of the other parameters. The class of split permutation graphs is $K_{2,2}$-free, has unbounded clique-width \cite{KLM14} but bounded twin-width \cite{BKTW21}, mim-width \cite{BV13} and hence sim-width. Therefore, clique-width does not become equivalent to twin-width, mim-width or sim-width. The class of chordal graphs is $K_{2,2}$-free, has unbounded mim-width but sim-width at most~$1$~\cite{KKST17}. Therefore, mim-width does not become equivalent to sim-width. 
\end{proof}

Note that the relationship between sim-width and tree-independence number is not covered by \Cref{comp} and indeed corresponds to \Cref{kttfreeopen}.
 
Let us now analyse the incomparable pairs and check whether they become comparable or equivalent, starting with pairs involving twin-width.

\begin{lemma}\label{incomp} Even when restricted to $K_{2,2}$-free graphs, twin-width is incomparable with
  \begin{itemize}
    \item mim-width,
    \item sim-width, and
    \item tree-independence number.
  \end{itemize}
\end{lemma}

\begin{proof} The class of walls is $K_{2,2}$-free, is not $(\tw,\omega)$-bounded \cite{DMS21a}, and hence has unbounded tree-independence number \cite{DMS24}, unbounded mim-width \cite{BHMPP22}, and unbounded sim-width by \cref{simwall}, but has bounded twin-width \cite{BKTW21}. On the other hand, the class of chordal graphs is $K_{2,2}$-free, has unbounded twin-width \cite{BGKTW21}, but bounded tree-independence number \cite{DMS24}, and hence has bounded sim-width. Finally, the class of interval graphs is $K_{2,2}$-free, has unbounded twin-width \cite{BGKTW21}, but bounded mim-width~\cite{BV13}.
\end{proof}

We finally consider pairs involving tree-independence number and show that tree-independence number dominates mim-width on $K_{t,t}$-free graphs (\Cref{tinmim}). The following two lemmas will be used. 

\begin{lemma}\label{deggeq}
Let $j$ and $\ell$ be positive integers. Let $G$ be a graph and let $U$ and $V$ be disjoint subsets of $V(G)$ such that each $u \in U$ has at least one neighbour in $V$ while each $v \in V$ has at most $j$ neighbours in $U$. If $|U| \geq 2j\ell$, then there exist $X \subseteq U$ and $Y \subseteq V$ such that $|X|=|Y| = \ell$ and $G[X,Y] \cong \ell P_2$.
\end{lemma}

\begin{proof} We proceed by induction on $\ell$. The base case $\ell=1$ is trivial. Let $\ell' > 1$ and suppose that the statement holds for each $\ell < \ell'$. We show it holds for $\ell'$. Therefore, let $|U| \geq 2j\ell'$. Pick $x \in U$ such that $|N_V(x)|$ is minimum and let $y \in N_V(x)$. Then, at most $j-1$ vertices of $U \setminus \{x\}$ have the same neighbourhood in $V$ as $x$, or else $|N_U(y)|\geq j+1$, a contradiction. Let $U'' = \{u \in U \setminus \{x\} :  N_{V}(u) \neq N_V(x)\}$. For every vertex $u' \in U''$, we have that $N_V(u') \setminus N_V(x) \neq \varnothing$ by minimality of $|N_V(x)|$. Note that $|U''| \geq |U\setminus\{x\}| - (j-1) = |U| - j \geq 2j\ell'-j$. Let $U' = U'' \setminus N_U(y)$. As $|N_U(y)|\leq j$, we have that $|U'| \geq 2j\ell'-j - j = 2j(\ell'-1)$. Let $V' = V \setminus N_V(x)$. Consider now the graph $G[U' \cup V']$. Each vertex in $U'\subseteq U''$ has at least one neighbour in $V'$ and each vertex in $V'$ has at most $j$ neighbours in $U'$. Therefore, by the induction hypothesis, there exist $X' \subseteq U'$ and $Y' \subseteq V'$ such that $G[X',Y'] \cong (\ell'-1)P_2$. Note further that $y$ is anticomplete to $U' = U''  \setminus N_U(y)$ and $x$ is anticomplete to $V' = V  \setminus N_V(x)$, so taking $X= X' \cup \{x\}$ and $Y = Y' \cup \{y\}$ completes the proof.
\end{proof}

\begin{lemma}\label{fin} Let $m$ be a fixed positive integer. For each positive integer $n$ and $k$, let $f(n,k) = 2^{n+k}$ and $g_m(n,k) = mk^{n}$.
  Let $G$ be a graph and let $U$ and $V$ be disjoint subsets of $V(G)$, where $U$ is an independent set. Suppose, for some positive integers $n$ and $k$, that $|U| \geq f(n,k)$ and that $\alpha(G[N_V(u)]) \geq g_m(n,k)$ for each $u \in U$. Then, one of the following occurs:
\begin{enumerate}
\item there exist $X \subseteq U$ and $Y \subseteq V$ such that $|X|=|Y| = k$ and $G[X,Y] \cong kP_2$; or
\item there exist independent sets $X \subseteq U$ and $Y \subseteq V$ such that $|X|=n,|Y| = m$ and $G[X\cup Y] \cong K_{n,m}$.
\end{enumerate}
\end{lemma}

\begin{proof} We proceed by double induction on $n$ and $k$. The base case $n = 1$ or $k=1$ is an easy exercise left to the reader. Suppose that the statement is true for $(n,k)=(n',k'-1)$, with $n' \geq 1$ and $k' > 1$, and for $(n,k)=(n'-1,k')$, with $n' > 1$ and $k' \geq 1$. We show that the statement is true for $(n,k)=(n',k')$. Therefore, let $|U| \geq f(n',k')$ and assume $\alpha(G[N_V(u)]) \geq g_m(n',k')$ for each $u \in U$. 

  Pick $x \in U$ and let $V' = N_V(x)$. Suppose first that there exist at least $f(n'-1,k')$ vertices $u' \in U\setminus \{x\}$ such that $\alpha(G[N_{V'}(u')]) = \alpha(G[N_V(u') \cap V']) \geq g_m(n'-1,k')$. Let $U'$ be the set of such vertices. Hence, $|U'| \geq f(n'-1,k')$. By the induction hypothesis for $(n,k)=(n'-1,k')$, either there exist $X' \subseteq U'$ and $Y' \subseteq V'$ such that $|X'|=|Y'| = k'$ and $G[X',Y'] \cong k'P_2$, or there exist independent sets $X' \subseteq U'$ and $Y' \subseteq V'$ such that $|X'|=n'-1$, $|Y'| = m$ and $G[X'\cup Y'] \cong K_{n'-1,m}$. If the former occurs, set $X=X'$ and $Y=Y'$. If the latter occurs, set $X=X'\cup \{x\}$ and $Y =Y'$. In either case, it is easy to see that the statement holds for $(n,k)=(n',k')$ with the chosen $X$ and $Y$. 

Suppose instead that fewer than $f(n'-1,k')$ vertices $u' \in U\setminus \{x\}$ satisfy $\alpha(G[N_V(u') \cap V']) \geq g_m(n'-1,k')$. This implies that the number of vertices $a \in U\setminus\{x\}$ such that $\alpha(G[N_V(a) \cap V']) < g_m(n'-1,k')$ is at least $(f(n',k') - 1) - (f(n'-1,k') - 1) = 2^{n'+k'}-2^{n'-1+k'} = 2^{n'+k'-1} = f(n',k'-1)$. Let $A$ be the set of such vertices. Hence, $|A| \geq f(n',k'-1)$. Let now $B = V \setminus V' = V \setminus N_V(x)$. Observe that, for each $a \in A$, we have
\[\begin{array}{lcl}
\alpha(G[N_B(a)]) &= &\alpha(G[N_V(a)\setminus V']) \\[1pt]
& \geq & \alpha(G[N_V(a)]) - \alpha(G[N_V(a) \cap V'])\\[1pt]
& > &g_m(n',k') - g_m(n'-1,k')\\[1pt]
& = & m(k')^{n'-1}(k'-1)\\[1pt]
& > & m(k'-1)^{n'-1}(k'-1)\\[1pt]
& = & g_m(n',k'-1).
\end{array}\]
By the induction hypothesis for $(n,k)=(n',k'-1)$, either there exist $X' \subseteq A$ and $Y' \subseteq B$ such that $|X'|= |Y'| = k'-1$ and $G[X',Y'] \cong (k'-1)P_2$, or there exist independent sets $X' \subseteq A$ and $Y' \subseteq B$ such that $|X'|=n'$, $|Y'| = m$ and $G[X'\cup Y'] \cong K_{n',m}$. If the latter occurs, set $X=X'$ and $Y=Y'$ and the statement holds for $(n,k)=(n',k')$ with the chosen $X$ and $Y$. If however the former occurs, first set $X=X'\cup \{x\}$. Observe now that, since $x$ is anticomplete to $Y' \subseteq B =V \setminus V'$, it is sufficient to find a neighbour of $x$ in $V$ which is anticomplete to $X'$. Since $x\in U$, we have that $\alpha(G[V']) \geq g_m(n',k')$. Let $I$ be an independent set of $G[V']$ of size at least $g_m(n',k')$. By definition of $A$, each $x' \in X' \subseteq A$ is such that $\alpha(G[N_V(x') \cap V']) < g_m(n'-1,k')$ and so, for each $x' \in X'$, we have $|N_V(x') \cap I| < g_m(n'-1,k')$. Therefore, since $|X'| = k' -1$, we have that $|N_V(X') \cap I| < k'g_m(n'-1,k')$. However, $|I| \geq g_m(n',k') = k'g_m(n'-1,k')$, and so $I \setminus N_V(X') \neq \varnothing$. This implies that at least one vertex in $V'$ is anticomplete to $X'$. Pick such a vertex $y$ and set $Y = Y'\cup \{y\}$. The statement then holds for $(n,k)=(n',k')$ with the chosen $X$ and $Y$.
\end{proof}

We are finally ready to prove \Cref{tinmim}, which we restate for convenience. 

\tinmim*

\begin{proof} Let $f(n,k) = 2^{n+k}$ and $g_m(n,k) = mk^{n}$. We begin by describing an algorithm that constructs a tree decomposition $\mathcal{T}$ of $G$ from the branch decomposition $(T,\delta)$ of $G$. The bags of $\mathcal{T}$ will be constructed recursively. We first need to introduce some notation (see also \Cref{f-treedec}). 

\begin{figure}
\begin{center}
\includegraphics[scale=0.9]{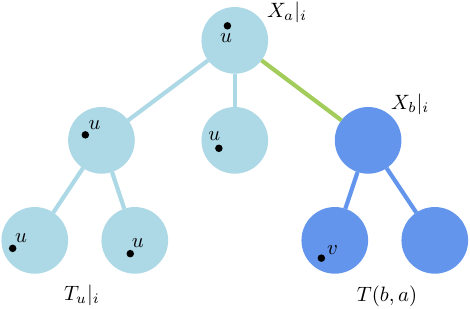}
\caption{Construction of a tree decomposition $\mathcal{T}$ of $G$ from a branch decomposition $(T, \delta)$ of $G$. The edge $ab$ of $T$ (in green) touches $T_u|_i$ (in light blue) at $a$. Informally, the set $N(u, b, a)|_i$ consists of all neighbours $v$ of $u$ in~$G$, hosted by $T(b, a)$, and such that the edge $uv$ does not yet satisfy (T2).}
\label{f-treedec}
\end{center}
\end{figure}

For each $t\in V(T)$ and $i \in \mathbb{N}$, we let $X_t|_i$ denote a particular subset of $V(G)$, where $i$ represents the step of the recursion (think of $X_t|_i$ as a bag assigned to $t$ at step $i$). Given a pair $(T,\{X_t|_i\}_{t \in V(T)})$ satisfying (T1) and (T3) in the definition of a tree decomposition (but not necessarily (T2)) and $v \in V(G)$, we denote by $T_v|_i$ the subtree of $T$ induced by the set $\{t\in V(T) : v \in X_t|_i\}$. Let $e = ab$ be an edge of $T$. Deleting $e$ splits $T$ into two subtrees, $T(a,b)$ and $T(b,a)$, where $a \in V(T(a,b))$ and $b \in V(T(b,a))$. We say that $T(a,b)$ (resp. $T(b,a)$) \textit{hosts} $v \in V(G)$ if $\delta(v)$ is a leaf of $T(a,b)$ (resp. $T(b,a)$). For $u \in V(G)$, the edge $ab \in E(T)$ \textit{touches $T_u|_i$ at $a$} if $a \in V(T_u|_i)$ and $b \not \in V(T_u|_i)$. If the edge $ab\in E(T)$ touches $T_u|_i$ at $a$, we let \[N(u,b,a)|_i = \{v \in N_G(u)\setminus X_a|_i : T(b,a) \ \mbox{hosts} \ v\}.\]

\noindent
\textbf{Algorithm.} 
We now describe the algorithm. We first pre-process $T$ by recursively contracting edges having at least one endpoint of degree $2$.
So we may assume that all internal nodes of $T$ have degree $3$.
We let $i$ represent a step counter. 

\medskip
\noindent
Set $i=0$. Set $X_t|_0 = \{\delta^{-1}(t)\}$ if $t$ is a leaf of $T$, and set $X_t|_0 = \varnothing$ if $t$ is an internal node of $T$. For each triple $(u,b,a)$, where $u \in V(G)$ and $ab$ touches $T_u|_0$ at $a$, compute $N(u,b,a)|_0$.

\medskip
\noindent
While there exists a triple $(u,b,a)$, where $u \in V(G)$, the edge $ab$ touches $T_u|_i$ at $a$, and $N(u,b,a)|_i$ contains an independent set of $G$ of size $g_m(n,k)$, do:
\begin{adjustwidth}{40pt}{40pt}
Pick an arbitrary such triple $(u,b,a)$. Set $X_b|_{i+1} = X_b|_{i} \cup \{u\}$, label $u$ a \textit{bad vertex with respect to $X_b$} and, for each $c\neq b$, set $X_c|_{i+1} = X_c|_{i} $. Compute $N(u,b,a)|_{i+1}$ for each triple $(u,b,a)$ where $u \in V(G)$ and $ab$ touches $T_u|_{i+1}$ at $a$. Set $i = i+1$. 
\end{adjustwidth}

\medskip
\noindent
While there exists a triple $(u,b,a)$, where $u \in V(G)$, the edge $ab$ touches $T_u|_i$ at $a$, and $N(u,b,a)|_i \neq \varnothing$, do: 
\begin{adjustwidth}{40pt}{40pt}
Pick an arbitrary such triple $(u,b,a)$. Set $X_b|_{i+1} = X_b|_{i} \cup \{u\}$, label $u$ a \textit{good vertex with respect to $X_b$} and, for each $c\neq b$, set $X_c|_{i+1} = X_c|_{i} $. Compute $N(u,b,a)|_{i+1}$ for each triple $(u,b,a)$ where $u \in V(G)$ and $ab$ touches $T_u|_{i+1}$ at $a$. Set $i = i+1$.
\end{adjustwidth}

\medskip
\noindent
Return $\mathcal{T} = (T,\{X_t|_i\}_{t \in V(T)})$. 

\begin{observation} Observe that no vertex added to a bag in the second loop is bad. Indeed, suppose that the first loop stops at step $i$. Then, for any triple $(u,b,a)$ where $u \in V (G)$ and the edge $ab$ touches $T_u|_i$ at $a$, we have that $N(u,b,a)|_i$ does not contain an independent set of $G$ of size $g_m(n,k)$. Let now $i' \geq i$ and fix an arbitrary triple $(u,b',a')$ such that $u \in V(G)$ and the edge $a'b'$ touches $T_{u}|_{i'}$ at $a'$. Consider the edge $ab$ of $T$ that touches $T_{u}|_{i}$ at $a$ and which is closest to $a'b'$. It holds that $N(u,b',a')|_{i'} \subseteq N(u,b,a)|_i$, and so $N(u,b',a')|_{i'}$ does not contain an independent set of $G$ of size $g_m(n,k)$.
\end{observation}

\noindent
\textbf{Running time analysis.} We now analyse the running time of the algorithm. Let $z = |V(G)|$ so, by definition, $z$ is the number of leaves of $T$. Since each vertex of $T$ has degree either $1$ or $3$, the number of internal vertices of $T$ is $z-2$ and $T$ has at most $2z - 3$ edges. Therefore, at each iteration of the while-do loops, there are $O(z^2)$ triples to be checked. For each such triple $(u, b, a)$, the set $N(u,b,a)|_i$ has size at most $z$ and so checking whether it contains an independent set of $G$ of size $g_m(n,k)$ can be done in $O(z^{g_m(n,k)})$ time. Since computing all the sets $N(u,b,a)|_{i+1}$ takes $O(z^{3})$ time, each iteration of the while-do loops can be done in $O(z^{g_m(n,k)+2} + z^{3}) = O(z^{g_m(n,k)+2})$ time. Observe now that, after each iteration of either loop, if $(u, b, a)$ is the chosen triple, then $T_u|_i$ is extended into $T_u|_{i+1}$ by the addition of one node of $T$. For each $u \in V(G)$, at most $O(z)$ such extensions are possible and so the total number of iterations of the while-do loops is $O(z^2)$. Therefore, the running time of the algorithm is $O(z^{g_m(n,k)+4})$.

\medskip
\noindent
\textbf{Correctness.} We first show that the algorithm indeed outputs a tree decomposition of $G$. Suppose that we stop at step~$s$, so $N(u,b,a)|_s = \varnothing$ for all triples $(u,b,a)$ such that $u \in V(G)$ and $ab$ touches $T_u|_s$ at $a$. We show that $\mathcal{T} = (T,\{X_t|_s\}_{t \in V(T)})$ is indeed a tree decomposition of $G$; namely, it satisfies (T1), (T2) and (T3). Since $\delta$ is a bijection from $V(G)$ to the leaves of $T$ and since $\delta^{-1}(t) \in X_t|_s$ for each leaf $t$ of $T$, (T1) holds. Consider now (T3). Let $u$ be an arbitrary vertex of $G$. We claim that $T_u|_s$ is connected. Observe that, at each step $i+1$ (for $i \geq 0$), we add $u$ to a bag $X_b|_i$ only if there exists an edge $ab \in E(T)$ that touches $T_u|_i$ at $a$. This implies that $T_u|_{i+1}$ is obtained from $T_u|_{i}$ by adding the node $b$ which is adjacent to $a \in T_u|_i$. Hence, if $T_u|_{i}$ is connected, then $T_u|_{i+1}$ is connected. Since $T_u|_{0}$ is connected, the same holds for $T_u|_{s}$. Finally, consider (T2). Suppose (T2) does not hold, so there exists $uv \in E(G)$ such that no bag in $\{X_t|_s\}_{t \in V(T)}$ contains both $u$ and $v$. Then $T_u|_s$ and $T_v|_s$ share no common nodes. Since $T$ is connected and, by (T3), $T_u|_s$ and $T_v|_s$ are connected subgraphs of $T$, there is a path in $T$ from $T_u|_s$ to $T_v|_s$. The first edge of this path when traversing from $T_u|_s$ to $T_v|_s$, say $ab$, touches $T_u|_s$ at $a$. Then $v \in N(u,b,a)|_s$, so $N(u,b,a)|_s \neq \varnothing$, contradicting that the algorithm terminates at step $s$.

Finally, we show that $\alpha(\mathcal{T}) < 3f(n,k)+6g_m(n+1,k)$. Suppose, to the contrary, that $\alpha(\mathcal{T}) \geq 3f(n,k)+6g_m(n+1,k)$. Then, there exists a bag $X_t \subseteq V(G)$, with $t \in V(T)$, such that $\alpha(G[X_t]) \geq 3f(n,k)+6g_m(n+1,k)$. Let $P \subseteq X_t$ be an independent set of $G$ such that $|P| \geq 3f(n,k)+6g_m(n+1,k)$. We claim that there exists $t' \in V(T)$ with $t't \in E(T)$ and such that $T(t',t)$ hosts at least $|P|/3$ vertices of~$P$. If $t$ is a leaf of $T$, it is enough to take as $t'$ the unique neighbour of $t$ in $T$, as $T(t',t)$ hosts $|P| - 1 \geq |P|/3$ vertices of $P$. If $t$ is an internal node of $T$, then $T-t$ is the disjoint union of three trees (recall that after our pre-processing all internal nodes have degree $3$). One of these trees, say $T_1$, must host at least $|P|/3$ vertices of $P$, and we let $t' \in V(T_1)$ be the unique vertex such that $t't \in E(T)$. This shows our claim. Let now $P' \subseteq P$ be the subset of vertices of $P$ hosted by $T(t', t)$. Hence, $|P'| \geq |P|/3 = f(n,k)+2g_m(n+1,k)$. Since each vertex of $P' \subseteq X_t$ has been labelled either good or bad with respect to $X_t$, we have that either at least $f(n,k)$ vertices of $P'$ are bad, or at least $2g_m(n+1,k)$ vertices of $P'$ are good.
 
Suppose first that at least $f(n,k)$ vertices of $P'$ are bad. Let $U \subseteq P' \subseteq P$ be the set of such vertices and let $V \subseteq V(G)$ be the set of vertices not hosted by $T(t',t)$. Pick an arbitrary $u\in U$ and suppose that $u$ has been added to $X_t$ at step $i+1$ for some $i\geq 0$. Then, $\alpha(G[N(u,t,t')|_{i}]) \geq g_m(n,k)$, where $N(u,t,t')|_{i} \subseteq N_V(u)$. Therefore, we found an independent set $U$ of $G$ disjoint from $V$ and such that $|U| \geq f(n,k)$ and, for each $u \in U$, we have $\alpha(G[N_V(u)]) \geq \alpha(G[N(u,t,t')|_{i}]) \geq g_m(n,k)$. Since $G$ is $K_{n,m}$-free, \Cref{fin} implies that there exist $X \subseteq U$ and $Y \subseteq V$ such that $|X|=|Y| = k$ and $G[X,Y] \cong kP_2$. But then, letting $e = tt'$, we obtain that $\cutmim_{G}(A_e,\overline{A_e}) \geq k$, contradicting that $\mimw_{G}(T,\delta) < k$.

Suppose instead that at least $2g_m(n+1,k)$ vertices of $P'$ are good. Let $U = \{u_1, u_2, \dots\}\subseteq P' \subseteq P$ be the set of such vertices. Hence, $|U| \geq 2g_m(n+1,k)$. For each $u_j \in U$, let $i_j+1$ be the step $u_j$ has been added to $X_t$, for some $i_j \geq 0$. Let $u_1$ be the first vertex of $U$ added to $X_t$, at step $i_{1}+1$. For $u_1$ to be added to $X_t$ at step $i_1+1$ as a good vertex, it must be that $N(u_1,t,t')|_{i_1} \neq \varnothing$ but $\alpha(G[N(u_1,t,t')|_{i_1}]) < g_m(n,k)$, for the edge $tt' \in E(T)$ touching $T_{u_1}|_{i_1}$ at $t'$. Moreover, for every other triple $(y,v,v')$, where $y \in V(G)$ and $vv'$ touches $T_{y}|_{i_1}$ at $v'$, it must be that $\alpha(G[N(y,v,v')|_{i_1}]) < g_m(n,k)$. 

Now let $V = \bigcup_{u_j \in U}N(u_j,t,t')|_{i_j}$. We claim that each $v \in V$ is adjacent to fewer than $g_m(n,k)$ vertices of $U$. Suppose, to the contrary, that there exists $v \in N(u_j,t,t')|_{i_j}$, for some $u_j \in U$, which is complete to a subset $U' \subseteq U$ with $|U'| \geq g_m(n,k)$. Then, as $v \in N(u_j,t,t')|_{i_j}$ for some $u_j \in U$, it follows that $v \not \in X_{t'}|_{i_j}$. Since $i_j \geq i_1$, this implies that $v \notin X_{t'}|_{i_1}$ and so $T_v|_{i_1}$ does not contain the node $t'$. Since $v \in N(u_j,t,t')|_{i_j}$, we have that $v$ is hosted by $T(t,t')$. Now let $ab$ be the first edge of the shortest path in $T$ from $T_v|_{i_1}$ to $t'$. Clearly, $ab$ touches $T_v|_{i_1}$, say without loss of generality at $a$. Note that $T(t',t)$ is a subtree of $T(b,a)$. Since $u_1$ is the first vertex of $U$ added to $X_t$, no other $u_j \in U$ is added to $X_t$ at step $i_1$. We now claim that no vertex $u_j \in U$ belongs to $X_a|_{i_1}$. Observe first that $u_j$ is hosted by $T(t',t)$, i.e. $\delta(u_j)$ is a leaf of $T(t',t)$, and $a$ is a node of $T(t,t')$. Since $T$ is a tree, there is a unique path in $T$ from $\delta(u_j)$ to $a$ and such path contains $t$. By definition of $u_1$, we have that $u_j\not \in X_t|_{i_1}$. But if $u_j \in X_a|_{i_1}$ (in particular, $a \neq t$), then $T_{u_j}|_{i_1}$ is not connected, a contradiction. Therefore, for each $u' \in U'$, it must be that $u' \in N_G(v)\setminus X_a|_{i_1}$. This implies that $U' \subseteq N(v,b,a)|_{i_1}$ and so $\alpha(G[N(v,b,a)|_{i_1}]) \geq \alpha(G[U']) = |U'| \geq g_m(n,k)$. But then, at step $i_1 + 1$, the vertex $v$ should have been added to $X_b|_{i_1}$ as a bad vertex instead of adding $u_1$ to $X_t|_{i_1}$ as a good vertex, a contradiction.

By the previous paragraph, each $v \in V$ has fewer than $g_m(n,k)$ neighbours in $U$. Moreover, for each $u_j \in U$, the vertex $u_j$ has been added to $X_t$ and so $u_j$ has a neighbour in $V$. Therefore, $U$ and $V$ satisfy the conditions of \Cref{deggeq} with $j = g_m(n,k)$ and $\ell=k$, and so there exist $X \subseteq U$ and $Y \subseteq V$ such that $|X|= |Y| = k$ and $G[X,Y] \cong kP_2$. Since $U \subseteq T(t',t)$ and $V \subseteq T(t,t')$, by letting $e = tt'$ we have that $\cutmim_{G}(A_e,\overline{A_e}) \geq k$, contradicting that $\mimw_{G}(T,\delta) < k$.

To summarize, $\alpha(\mathcal{T}) < 3f(n,k)+6g_m(n+1,k) = 6(2^{n+k-1} + mk^{n+1})$, as desired.  
\end{proof}

We are finally ready to prove \Cref{t-main2}, which we restate for convenience.

\tmaintwo*

\begin{proof} It follows from \Cref{comp,incomp}, \Cref{tinmim} and the following observation showing that tree-independence number is more powerful than mim-width for $K_{t,t}$-free graphs when $t \ge 2$: Chordal graphs are $K_{2,2}$-free, have tree-independence number $1$ \cite{DMS24} but unbounded mim-width \cite{KKST17}.
\end{proof}

\section{Width Parameters and Graph Powers: The Proof of \Cref{powersim}}\label{s-powers}

In this short section we prove \Cref{powersim} and observe that a result similar to \Cref{powersim} cannot hold for even powers.

\powersim*

\begin{proof} Clearly, $(T,\delta)$ is a branch decomposition of $G^r$. Suppose, to the contrary, that $\simw_{G^r}(T,\delta) > w$. Hence, $r \geq 3$. Let $e \in E(T)$ be such that $\cutsim_{G^r}(A_e,\overline{A_e}) \geq w+1$. There exist independent sets $X=\{x_1,\ldots,x_{w+1}\} \subseteq A_e$ and $Y=\{y_1,\ldots,y_{w+1}\} \subseteq \overline{A_e}$ of $G^r$ such that, for each $i,j \in \{1, \ldots, w+1\}$, $x_i$ is adjacent to $y_j$ if and only if $i=j$. Since $x_i$ is adjacent to $y_i$ in $G^r$ for each $i \in \{1, \ldots, w+1\}$, there exists a path $P_i = v(i,0)v(i,1)v(i,2)\cdots v(i,a_i)$ in $G$, with endpoints $v(i,0) = x_i$ and $v(i,a_i)=y_i$, for some $a_i \leq r$. We claim that, for each $i \neq j$, $V(P_i) \cap V(P_j) = \varnothing$ and no vertex of $P_i$ is adjacent to a vertex of $P_j$ in $G$. Suppose, to the contrary, that there exist $v(i,k) \in V(P_i)$ adjacent to $v(j,l) \in V(P_j)$ (the proof that $V(P_i) \cap V(P_j) = \varnothing$ is similar and left to the reader).

Suppose first that $k+l < r$. Then, there is a walk $v(i,0)v(i,1)\cdots v(i,k)v(j,l)v(j,l-1)\cdots v(j,0)$ in $G$ of length at most $r$ and hence a $v(i,0),v(j,0)$-path in $G$, from which $x_i=v(i,0)$ is adjacent to $x_j=v(j,0)$ in $G^r$, a contradiction. Suppose now that $k+l > r$. Then, there is a walk $v(i,a_i)v(i,a_i-1)\cdots v(i,k)v(j,l)v(j,l+1)\cdots v(j,a_j)$ in $G$ of length at most $r$, from which $y_i = v(i,a_i)$ is adjacent to $y_j = v(j,a_j)$ in $G^r$, a contradiction. Suppose finally that $k+l = r$. Since $r$ is odd, either $k < l$ or $k > l$. Without loss of generality, $k < l$. Then, there is a walk  $v(i,0)v(i,1)\cdots v(i,k)v(j,l)v(j,l+1)\cdots v(j,a_j)$ in $G$ of length at most $r$, and so $x_i = v(i,0)$ is adjacent to $y_j = v(j,a_j)$ in $G^r$, a contradiction. 

Consider now, for each $i \in \{1, \ldots, w+1\}$, the path $P_i = v(i,0)v(i,1)v(i,2)\cdots v(i,a_i)$. Since $x_i = v(i,0) \in A_e$ and $y_i = v(i,a_i) \in \overline{A_e}$, there exists an integer $b_i \leq a_i-1$ such that $v(i,b_i) \in A_e$ and $v(i,b_i+1) \in \overline{A_e}$. Let $x'_i = v(i,b_i)$ and $y'_i = v(i,b_i+1)$. Clearly, each $x'_i$ is distinct from and adjacent to $y'_i$. Moreover, by the paragraph above, each $x'_i$ is distinct from and non-adjacent to both $x'_j$ and $y'_j$, for $i \neq j$. We then have that $X' = \{x'_1,\ldots,x'_{w+1}\} \subseteq A_e$ and $Y' = \{y'_1,\ldots,y'_{w+1}\} \subseteq \overline{A_e}$ are independent sets of $G$ and $G[X',Y'] \cong (w+1)P_2$. Therefore, $\cutsim_{G}(A_e,\overline{A_e}) \geq w+1$, a contradiction.
\end{proof}

Lima et al.~\cite{LMM24} showed that, for every fixed even integer $r \geq 2$ and for every graph $H$, there exists a chordal graph $G$ such that $G^r$ contains an induced subgraph isomorphic to $H$. Since walls have arbitrarily large sim-width, their result immediately implies the following:

\begin{proposition}\label{p-new} 
For every even integer $r \geq 2$ and every integer $w \geq 1$, there exists a graph $G$ such that $\simw(G) = 1$ while $\simw(G^r) \geq w$. In particular, for every fixed even integer $r \geq 2$, there is no function $f$ such that $\simw(G^r) \leq f(\simw(G))$ for all graphs $G$. 
\end{proposition}

\section{Concluding Remarks and Open Problems}\label{s-con}

In \Cref{t-main,t-main3,t-main2}, we investigated the relationships between six width parameters (treewidth, clique-width, twin-width, mim-width, sim-width and tree-independence number) when restricted to $K_{t,t}$-subgraph-free graphs, line graphs and $K_{t,t}$-free graphs in order to examine to what extent relationships between non-equivalent width parameters may change (see \Cref{fig:gadham}). In this way we also extended and generalized several known results from the literature. Moreover, in the case that two parameters become comparable or equivalent on one of these graph classes, we showed how to obtain computable functions witnessing this. 

Arguably, the main unresolved problem is the (only) missing case in Figure~\ref{fig:gadham}, which corresponds to the following question already stated in \Cref{sec:intro}.

\kttfreeopen*

We first observe a consequence of a positive answer to \Cref{kttfreeopen}. If tree-independence number dominates sim-width for the class of $K_{t,t}$-free graphs then, in order to prove \Cref{twomegaconj}, it suffices to show that every hereditary $(\tw, \omega)$-bounded graph class has bounded sim-width. The question of whether every $(\tw, \omega)$-bounded graph class has bounded sim-width was first asked in \cite{MY23}. Note also that complete bipartite graphs have sim-width $1$ but are not $(\tw, \omega)$-bounded. 

We now turn to possible algorithmic consequences of \Cref{tinmim}, and of its extension contingent to a positive answer to \Cref{kttfreeopen}. Computing optimal decompositions for a certain width is in general an $\mathsf{NP}$-hard problem (see \cite{SV16} and references therein). However, in some cases, there exist exact or approximation algorithms running in $\mathsf{FPT}$ or $\mathsf{XP}$ time parameterized by the target width. This is exemplified by treewidth, which admits an exact $\mathsf{FPT}$ algorithm \cite{Bod96}. Most of the time, it is in fact sufficient to simply obtain an approximate decomposition: we seek an algorithm that, given a graph of width at most~$k$, outputs a decomposition of width at most $f(k)$, for some computable function $f$. 

It is known that rank-width admits an $\mathsf{FPT}$-approximation algorithm~\cite{OS06}, with the current best-known result, in terms of running time, being the following. 
Korhonen and Soko\l{}owski~\cite{KS24} showed that, for fixed~$k$, there exists an algorithm that, given an $n$-vertex $m$-edge graph $G$, in time $O_k(n^{1+o(1)}) + O(m)$, either decides that $\rw(G) > k$, or outputs a branch decomposition of $G$ of $\mathrm{cutrk}_{G}$-width at most $k$. Recently, it was shown that tree-independence number admits an $\mathsf{XP}$-approximation algorithm: Dallard et al.~\cite{DFGKM24} showed that, for fixed $k$, there exists an algorithm that, given an $n$-vertex graph $G$, in time $2^{O(k^2)}n^{O(k)}$, either decides that $\tin(G) > k$, or outputs a tree decomposition of $G$ with independence number at most $8k$. For other width parameters, such as mim-width and sim-width, it is a well-known open problem to obtain $\mathsf{XP}$-approximation algorithms.

\begin{open}[see, e.g., \cite{Jaf20,KKST17}]\label{XPsimwidth} Does there exist a computable function $f$ and an algorithm $A$ that, for fixed $k$ and given a graph $G$, in $\mathsf{XP}$ time parameterized by $k$, either decides that $\mimw(G) > k$ (or $\simw(G) > k$), or outputs a branch decomposition of $G$ of mim-width (or sim-width) at most $f(k)$?
\end{open}

Therefore, in contrast to algorithms on classes of bounded treewidth, clique-width or tree-independence number, algorithms for graph problems restricted to classes of bounded mim-width (or sim-width) require a branch decomposition of constant mim-width (or sim-width) as part of the input. Obtaining such branch decompositions in polynomial time has been shown possible for several special graph classes (see, for example, \cite{BV13,BHMPP22,KKST17}). 

One may also consider the problem of finding exact and approximation algorithms for computing optimal decompositions for a certain width parameterized by a parameter other than the target width. For example, Bodlaender and Kloks~\cite{BK96} obtained an $\mathsf{XP}$ algorithm for computing pathwidth when parameterized by the treewidth of the input graph, and it is not known whether this can be improved to $\mathsf{FPT}$ (see, e.g., \cite{GJNW22}). Eiben et al.~\cite{EGHJ22} showed that mim-width admits an exact $\mathsf{FPT}$ algorithm parameterized by the treewidth and the maximum degree of the input graph, and an exact $\mathsf{FPT}$ algorithm parameterized by the treedepth of the input graph. Groenland et al.~\cite{GJNW22} obtained a polynomial-time algorithm that approximates pathwidth to within a factor of $O(\tw(G)\sqrt{\log\tw(G)})$. Their key observation is that every graph with large pathwidth either has large treewidth or contains a subdivision of a large complete binary tree. This shows how the study of exact and approximation algorithms for computing optimal decompositions for a certain width is related to the study of obstructions to small width.

A straightforward consequence of \Cref{tinmim} is the following. Here the \textit{induced biclique number} of a graph $G$ is the largest $t \in \mathbb{N}$ such that $G$ contains $K_{t,t}$ as an induced subgraph. 

\begin{corollary}\label{paramdec} There exists an $\mathsf{XP}$-approximation algorithm for tree-independence number parameterized by rank-width and induced biclique number.
\end{corollary}

\begin{proof}
Let $G$ be the input graph, and let $t$ be the induced biclique number of $G$. We compute a branch decomposition of $G$ of $\mathrm{cutrk}_{G}$-width at most $2\rw(G)$ in time $2^{2^{O(\rw(G))}}n^2$ \cite{FK22}. This is also a branch decomposition of $G$ of mim-width less than $2\rw(G)+1$ (see, e.g., \cite[Lemma~2.4]{BK19}). We then run the algorithm from \Cref{tinmim} with this branch decomposition in input. It outputs, in $n^{O(t\rw(G)^{t})}$ time, a tree decomposition of $G$ with independence number $O(2^{\rw(G)} + t(2\rw(G)+1)^{t+1})$. 
\end{proof}

It is not immediately clear whether \Cref{paramdec} gives a conditional improved running time compared to the $2^{O(k^2)}n^{O(k)}$ $\mathsf{XP}$-approximation algorithm parameterized by tree-independence number \cite{DFGKM24} mentioned above. So we pose the following open problem.

\begin{open} For a graph $G$ with induced biclique number $t$ and rank-width $\rw(G)$, find an asymptotically tight upper bound on $\tin(G)$ in terms of $t$ and $\rw(G)$.
\end{open}

In a similar vein, we also observe an immediate consequence of a positive answer to Open Problems~\ref{kttfreeopen} and~\ref{XPsimwidth}. Suppose that \Cref{XPsimwidth} has a positive answer for sim-width. That is, suppose that there exists an $\mathsf{XP}$-approximation algorithm for sim-width parameterized by the target width. If, in addition, an algorithmic version of \Cref{kttfreeopen} has a positive answer (i.e., suppose that, given a $K_{t,t}$-free graph $G$ and a branch decomposition of $G$ with sim-width at most $k$, it is possible to compute a tree decomposition of $G$ with independence number at most $g(t,k)$, for some computable function $g$, in $\mathsf{XP}$ time parameterized by $k$ and $t$), then we can obtain an $\mathsf{XP}$-approximation algorithm for tree-independence number. Indeed, given an input graph $G$ and an integer $k$, we simply check whether $G$ is $K_{k+1,k+1}$-free. If not, then $\tin(G) > k$. Otherwise, $G$ is $K_{k+1,k+1}$-free. We now run the algorithm $A$ from \Cref{XPsimwidth}. Algorithm $A$ either decides that $\simw(G) > k$, and so $\tin(G) > k$, or outputs a branch decomposition of $G$ of sim-width at most $f(k)$, from which we build a tree decomposition of $G$ with independence number at most $g(k)$ in $\mathsf{XP}$ time parameterized by $k$. This would thus provide a different proof of the main result in \cite{DFGKM24}, although one should expect worse running time and approximation factor.

Other natural problems related to \Cref{t-main,t-main3,t-main2} consist of optimizing the bounding functions obtained therein, as those provided are most likely not optimal. A first attempt in this direction was made in \Cref{mimlinebw}, where we showed that, in fact, $\mimw(L(G))$ equals (up to a multiplicative constant) $\bw(G)$, for any graph $G$. It would be interesting to refine the bounds for the sim-width of a line graph, and a starting point is to improve \Cref{simline}.  

\begin{open}\label{openline}
Find an asymptotically optimal function $f$ such that $\simw(L(G)) \geq f(\tw(G))$ for any graph $G$. 
\end{open}

Any attempt to solve \Cref{openline} seems to need to avoid the use of the Grid-minor theorem. Another problem in this direction is whether it is possible to improve the bound in \Cref{tinmim}.

\begin{open} Is the exponential dependency of tree-independence number in mim-width and induced biclique number from \Cref{tinmim} necessary? 
\end{open}

In \Cref{simwclique}, we observed that $\left\lceil\frac{n}{6}\right\rceil \leq \simw(L(K_{n})) \leq \left\lceil\frac{2n}{3}\right\rceil$, for each $n > 12$. We ask to determine the exact value.

\begin{open} Determine the exact value of $\simw(L(K_{n}))$. 
\end{open} 

\noindent \textbf{Note.} Several results related to our paper have been announced after its submission for publication. First, \Cref{twomegaconj} was disproved by Chudnovsky and Trotignon~\cite{CT24}. Second, Abrishami et al.~\cite{ABC24} proved the following special case of \Cref{kttfreeopen} (incomparable with \Cref{tinmim}): every subclass of $K_{t,t}$-free graphs of bounded induced matching treewidth has bounded tree-independence number. Third, Dallard et al.~\cite{DKKM24} showed that, for every positive integer $n$, $\simw(L(K_{n})) \leq \left\lfloor\frac{n}{2}\right\rfloor$.

\section*{Acknowledgement}

The authors greatly appreciate the careful comments of the anonymous referees, in particular their suggestion of a simpler counterexample (\Cref{newcounter}) to \cite[Lemma~4.3.9]{Vat12}.

\bibliography{references}
\end{document}